\documentclass[amstex,12pt, amssymb]{article}

\usepackage{mathtext}
\usepackage[cp1251]{inputenc}
\usepackage[T2A]{fontenc}
\usepackage[dvips]{graphicx}
\usepackage{amsmath}
\usepackage{amssymb}
\usepackage{amsxtra}
\usepackage{latexsym}
\usepackage{ifthen}

\newcounter{lemma}[section]

\newcounter{corollary}[section]

\newcounter{remark}[section]

\newcounter{theorem}[section]

\newcounter{proposition}[section]

\newcounter{example}

\numberwithin{equation}{section}

\begin{document}

\markboth{V. DESYATKA, E.~SEVOST'YANOV}{\centerline{ON THE BOUNDARY
BEHAVIOR ...}}

\def\cc{\setcounter{equation}{0}
\setcounter{figure}{0}\setcounter{table}{0}}

\overfullrule=0pt

%\normalsize\large

\author{VICTORIA DESYATKA, EVGENY SEVOST'YANOV}

\title{
{\bf CARATH\'{E}ODORY BOUNDARY EXTENSIONS FOR GENERALIZED
QUASIREGULAR MAPPINGS}}

\date{\today}
\maketitle

%\large
\begin{abstract}
The manuscript is devoted to the boundary behavior of mappings with
bounded and finite distortion, which has been actively studied
recently. We consider mappings of domains of the Euclidean space
that satisfy the inverse Poletsky inequality with an integrable
majorant, are open, and discrete, and not necessarily preserve the
boundary of a domain. Under some conditions on the geometry of these
domains, it is proved that the specified mappings have a continuous
boundary extension. The result is valid even in a more general form,
when the majorant mentioned above is integrable over almost all
concentric spheres centered at each point. We also have proved some
results on equicontinuity of the family of above mappings in the
closure of a domain.
\end{abstract}


\section{Introduction}

One of the most important problems of analysis is a problem of a
continuous boundary extension of mappings. Research in this
direction was conducted by various authors, see e.g. \cite{A},
\cite{Car}, \cite{Cr}, \cite{GRSY}, \cite{MRSY$_2$},
\cite{Na$_1$}--\cite{Na$_2$}, \cite{Va$_1$} and  \cite{Vu}. In
particular, the following result on this topic is valid (see
\cite[Theorem~17.13]{Va$_2$}).

\medskip
{\bf Theorem~{\bf A.}}{\it\, Suppose that $f:D\rightarrow
D^{\,\prime}$ is a quasiconformal mapping and that $D$ has property
$P_1$ at $\in
\partial D.$ Then $C(f , b)$ contains at most one point at which
$D^{\,\prime}$ is finitely connected.}

\medskip
Here we say that {\it a domain $D$ has property $P_1$ at $b$} if the
following condition is satisfied: If $E$ and $F$ are connected
subsets of $D$ such that $b\in \overline{E}\cup \overline{F},$ then
$M(\Gamma(E, F, D))=\infty.$ In addition, $M(\Gamma(E, F, D))$
denotes the modulus of families of paths joining the sets $E$ and
$F$ in $D$ (see the definition below). Besides that, $C(f , b)$
denotes the cluster set of $f$ at the point~$b.$

\medskip An analogue of Theorem~{\bf A} was also obtained for quasiregular
mappings (which, in particular, are not homeomorphisms; see, for
example, \cite[Theorem~4.2]{Sr}).

\medskip
{\bf Theorem~{\bf B.}}{\it\, Let $f:D\rightarrow {\Bbb R}^n$ be
quasiregular mapping with $C(f, \partial D)\subset \partial f(D).$
If $D$ is locally connected at a point $b\in \partial D$ and
$D^{\,\prime}=f(D)$ is qc accessible at some point $y\in C(f, b),$
then $C(f, b)=\{y\}.$}

\medskip
We should note that both theorems make one important topological
assumption on the mapping: in the first theorem it is assumed to be
quasiconformal and, therefore, is a homeomorphism, which in turn
implies the condition $C(f, \partial D)\subset \partial f(D)$ (see,
e.g., \cite[Proposition~13.5]{MRSY$_2$}). In the second theorem, the
condition $C(f, \partial D)\subset \partial f(D)$ is directly
present. The question arises: how important (necessary) is this
condition for both theorems to remain valid? The general answer to
this question is, in a sense, quite simple: as follows, say, from
the Sokhotsky-Сasorati-Weierstrass theorem for analytic functions,
they, generally speaking, do not have a continuous extension even to
an isolated boundary point, if we discard this condition. Of course,
more complex types of boundaries, as well as more general classes of
quasiregular mappings, do not simplify the situation. However, as we
will see below, the condition $C(f, \partial D)\subset \partial
f(D)$ can be significantly weakened.

Some analogues of Theorems~\textbf{A} and~\textbf{B} were obtained
by the second co-author for more general classes of mappings than
quasiregular ones, see, for example, ~\cite{SSD} and \cite{Sev$_1$}.
These results were obtained under the boundary-preserving condition
for mappings. The aim of the present manuscript is to obtain
analogues of the same theorems when the boundary-preserving
condition is not assumed. It seems to us that the assertions
formulated below are new even for quasiregular mappings. However,
these results were formulated for more general classes of mappings
satisfying weighted inverse moduli conditions of Poletskii type.
Most known classes of mappings satisfy such conditions, in
particular, mappings with finite length distortion, see, for
example, \cite[Theorem~8.5]{MRSY$_2$}

\medskip
Let us recall the definitions. Let $y_0\in {\Bbb R}^n,$
$0<r_1<r_2<\infty$ and
\begin{equation}\label{eq1**}
A=A(y_0, r_1,r_2)=\left\{ y\,\in\,{\Bbb R}^n:
r_1<|y-y_0|<r_2\right\}\,.\end{equation}
Given sets $E,$ $F\subset\overline{{\Bbb R}^n}$ and a domain
$D\subset {\Bbb R}^n$ we denote by $\Gamma(E,F,D)$ a family of all
paths $\gamma:[a,b]\rightarrow \overline{{\Bbb R}^n}$ such that
$\gamma(a)\in E,\gamma(b)\in\,F $ and $\gamma(t)\in D$ for $t \in
(a, b).$ If $f:D\rightarrow {\Bbb R}^n,$ $y_0\in f( D)$ and
$0<r_1<r_2<d_0=\sup\limits_{y\in f(D)}|y-y_0|,$ then by
$\Gamma_f(y_0, r_1, r_2)$ we denote the family of all paths $\gamma$
in $D$ such that $f(\gamma)\in \Gamma(S(y_0, r_1), S( y_0, r_2),
A(y_0,r_1,r_2)).$ Let $Q:{\Bbb R}^n\rightarrow [0, \infty]$ be a
Lebesgue measurable function. We say that {\it $f$ satisfies
Poletsky inverse inequality} at the point $y_0\in \overline{f(D)},$
if the relation
\begin{equation}\label{eq2*A}
M(\Gamma_f(y_0, r_1, r_2))\leqslant \int\limits_{A(y_0,r_1,r_2)\cap
f(D)} Q(y)\cdot \eta^n (|y-y_0|)\, dm(y)
\end{equation}
holds for any Lebesgue measurable function $\eta:
(r_1,r_2)\rightarrow [0,\infty ]$ such that
\begin{equation}\label{eqA2}
\int\limits_{r_1}^{r_2}\eta(r)\, dr\geqslant 1\,.
\end{equation}
Recall that, the Poletsky classical inequality is the inequality of
the type
\begin{equation}\label{eq5}
M(f(\Gamma))\leqslant K\cdot M(\Gamma)\,,
\end{equation}
where $K\geqslant 1$ is some constant, and $\Gamma$ is an arbitrary
family of paths in $D.$ It was obtained for quasiregular mappings by
E.~Poletsky for quasiregular mappings in the early 70s last century,
see~\cite[Theorem~1]{Pol}. If there exists $f^{\,-1}=g,$ then
instead of~(\ref{eq5}) we may write $M(\Gamma)\leqslant K\cdot
M(g(\Gamma)),$ and this inequality is similar to~(\ref{eq2*A}) for
$Q\equiv K.$ This circumstance justifies our proposed name of the
relation~(\ref{eq2*A}).

\begin{remark}\label{rem1}
Note that all quasiregular mappings $f:D\rightarrow {\Bbb R}^n$
satisfy the condition
\begin{equation}\label{eq22}
M(\Gamma_f(y_0, r_1, r_2))\leqslant \int\limits_{f(D)\cap
A(y_0,r_1,r_2)} K_O\cdot N(y, f, D)\cdot \eta^n (|y-y_0|)\, dm(y)
\end{equation}
at each point $y_0\in \overline{f(D)}\setminus\{\infty\}$ with some
constant $K_O=K_O(f)\geqslant 1$ and an arbitrary
Lebesgue-dimensional function $\eta: (r_1,r_2)\rightarrow
[0,\infty],$ which satisfies condition~(\ref{eqA2}). Indeed,
quasiregular mappings satisfy the condition
\begin{equation}\label{eq24}
M(\Gamma_f(y_0, r_1, r_2))\leqslant \int\limits_{f(D)\cap
A(y_0,r_1,r_2)} K_O\cdot N(y, f, D\setminus\{x_0\})\cdot
(\rho^{\,\prime})^n(y)\, dm(y)
\end{equation}
for an arbitrary function $\rho^{\,\prime}\in{\rm adm}\,
f(\Gamma_f(y_0, r_1, r_2)),$ see \cite[Remark~2.5.II]{Ri}. Put
$\rho^{\,\prime}(y):=\eta(|y-y_0|)$ for $y\in A(y_0,r_1,r_2)\cap
f(D),$ and $\rho^{\,\prime}(y)=0$ otherwise. By Luzin theorem, we
may assume that the function $\rho^{\,\prime}$ is Borel measurable
(see, e.g., \cite[Section~2.3.6]{Fe}). Due
to~\cite[Theorem~5.7]{Va$_2$},
$$\int\limits_{\gamma_*}\rho^{\,\prime}(y)\,|dy|\geqslant
\int\limits_{r_1}^{r_2}\eta(r)\,dr\geqslant 1$$
for each (locally rectifiable) path $\gamma_*$ in $\Gamma(S(y_0,
r_1), S(y_0, r_2), A(y_0, r_1, r_2)).$ By substituting the function
$\rho^{\,\prime}$ mentioned above into~(\ref{eq24}), we obtain the
desired ratio~(\ref{eq22}).
\end{remark}

\medskip
Recall that a mapping $f:D\rightarrow {\Bbb R}^n$ is called {\it
discrete} if the pre-image $\{f^{-1}\left(y\right)\}$ of each point
$y\,\in\,{\Bbb R}^n$ consists of isolated points, and {\it is open}
if the image of any open set $U\subset D$ is an open set in ${\Bbb
R}^n.$ Later, in the extended space $\overline{{{\Bbb R}}^n}={{\Bbb
R}}^n\cup\{\infty\}$ we use the {\it spherical (chordal) metric}
$h(x,y)=|\pi(x)-\pi(y)|,$ where $\pi$ is a stereographic projection
$\overline{{{\Bbb R}}^n}$ onto the sphere
$S^n(\frac{1}{2}e_{n+1},\frac{1}{2})$ in ${{\Bbb R}}^{n+1},$ namely,
\begin{equation}\label{eq3C}
h(x,\infty)=\frac{1}{\sqrt{1+{|x|}^2}}\,,\qquad
h(x,y)=\frac{|x-y|}{\sqrt{1+{|x|}^2} \sqrt{1+{|y|}^2}}\,, \quad x\ne
\infty\ne y
\end{equation}
(see \cite[Definition~12.1]{Va$_2$}). Further, the closure
$\overline{A}$ and the boundary $\partial A$ of the set $A\subset
\overline{{\Bbb R}^n}$ we understand relative to the chordal metric
$h$ in $\overline{{\Bbb R}^n}.$

\medskip
The boundary of $D$ is called {\it weakly flat} at the point $x_0\in
\partial D,$ if for every $P>0$ and for any neighborhood $U$
of the point $x_0$ there is a neighborhood $V\subset U$ of the same
point such that $M(\Gamma(E, F, D))>P$ for any continua $E, F\subset
D$ such that $E\cap\partial U\ne\varnothing\ne E\cap\partial V$ and
$F\cap\partial U\ne\varnothing\ne F\cap\partial V.$ The boundary of
$D$ is called weakly flat if the corresponding property holds at any
point of the boundary $D.$ Given a mapping $f:D\rightarrow {\Bbb
R}^n$, we denote
\begin{equation}\label{eq1_A_4} C(f, x):=\{y\in \overline{{\Bbb
R}^n}:\exists\,x_k\in D: x_k\rightarrow x, f(x_k) \rightarrow y,
k\rightarrow\infty\}
\end{equation}
and
\begin{equation}\label{eq1_A_5} C(f, \partial
D)=\bigcup\limits_{x\in \partial D}C(f, x)\,.
\end{equation}
In what follows, ${\rm Int\,}A$ denotes the set of inner points of
the set $A\subset \overline{{\Bbb R}^n}.$ Recall that the set
$U\subset\overline{{\Bbb R}^n}$ is neighborhood of the point $z_0,$
if $z_0\in {\rm Int\,}A.$

\medskip
Some analogues of the following result were established for the case
of homeomorphisms in \cite[Lemma~5.20, Corollary~5.23]{MRSY$_1$},
\cite[Lemma~6.1, Theorem~6.1]{RS} and \cite[Lemma~5, Theorem~3]{Sm}.
For open discrete and closed mappings, see, e.g., in \cite{SSD},
\cite{Sev$_1$} and \cite{Sev$_2$}.

\medskip
\begin{theorem}\label{th3}
{\it\, Let $D$ and $D^{\,\prime}$ be domains in ${\Bbb R}^n,$
$n\geqslant 2,$ and let $D$ be a domain with a weakly flat boundary.
Suppose that $f$ is open discrete mapping of $D$ onto $D^{\,\prime}$
satisfying the relation~(\ref{eq2*A}) at each point $y_0\in
\overline{D^{\,\prime}}.$ In addition, assume that the following
conditions are fulfilled:

\medskip
1) for each point $y_0\in \partial D^{\,\prime}\setminus\{\infty\}$
there is $0<r_0:=\sup\limits_{y\in D^{\,\prime}}|y-y_0|$ and any
$0<r_1<r_2<r_0:=\sup\limits_{y\in D^{\,\prime}}|y-y_0|$ there exists
a set $E_1\subset[r_1, r_2]$ of positive linear Lebesgue measure
such that $Q$ is integrable on $S(y_0, r)$ for $r\in E_1;$

\medskip
2) $C(f, \partial D)\subset E$ for some set $E\subset
\overline{D^{\,\prime}}$ which is closed in $\overline{{\Bbb R}^n},$
while $D^{\,\prime}$ is locally finitely connected with respect to
$E,$ in other words, for each point $z_0\in E$ and for any
neighborhood $U$ of this point there exists a neighborhood $V\subset
U$ of $z_0$ such that the set $V\cap (D^{\,\prime}\setminus E)$
consists of a finite number of components;

\medskip
3) the set $f^{\,-1}(E\cap D^{\,\prime})$ is nowhere dense in $D.$

Then the mapping $f$ has a continuous extension
$\overline{f}:\overline{D}\rightarrow\overline{D^{\,\prime}},$
moreover, $\overline{f}(\overline{D})=\overline{D^{\,\prime}}.$ }
\end{theorem}

\medskip
\begin{corollary}\label{cor1}
{\it\, The statement of Theorem~\ref{th3} remains true if, instead
of condition~1), a simpler condition holds: $Q\in
L^1(D^{\,\prime}).$}
\end{corollary}

\medskip
We also established a result on the equicontinuity of the mappings
discussed in Theorem~\ref{th3}. For quasiconformal mappings,
theorems of this kind were obtained, for example, in~\cite{NP}.
Below we give a version for mappings with unbounded characteristic.
As in Theorem~\ref{th3}, we do not assume that the mappings do not,
in general, preserve the boundary.

\medskip
Given a number $\delta>0,$ domains $D, D^{\,\prime}\subset{\Bbb
R}^n,$ $n\geqslant 2,$ a set $E\subset \overline{D^{\,\prime}},$
closed in $\overline{{\Bbb R}^n},$ and such that $\partial
D^{\,\prime}\subset E,$ a finite or countable set of points
$P=\{a_i\}_{i=1}^{\infty}\subset D^{\,\prime}\setminus E$ and a
Lebesgue measurable function $Q:D^{\,\prime}\rightarrow [0, \infty]$
denote by ${\frak S}^P_{E, \delta, Q}(D, D^{\,\prime})$ a family of
all open discrete mappings $f$ of $D$ onto $D^{\,\prime},$
satisfying the condition~(\ref{eq2*A}) for each $y_0\in
\overline{D^{\,\prime}}$ and such that

\medskip
1) $C(f, \partial D)\subset E,$

\medskip
2) the set $f^{\,-1}(E\cap D^{\,\prime})$ is nowhere dense in $D,$

\medskip
3) $h(f^{-1}(a_i), \partial D)\geqslant \delta$ for all $i\in {\Bbb
N}$ (here $h(A, B)$ denotes the chordal distance between sets in
$\overline{{\Bbb R}^n},$ $h(A, B)=\inf\limits_{x\in A, y\in B}h(x,
y)$).

The following theorem in some other versions (in particular, for
homeomorphisms and closed maps) was proved by the second co-author
in \cite{SSD}, \cite{SevSkv$_2$}, \cite{SevSkv$_1$} and
\cite{SevSkv$_3$}. The result formulated below takes into account
the most general situation, when the mapping is neither a
homeomorphism nor closed (=boundary preserving).

\medskip
\begin{theorem}\label{th2}
{\it\,Let $D$ and $D^{\,\prime}$ be domains in ${\Bbb R}^n,$
$n\geqslant 2,$ let $D$ be a domain with a weakly flat boundary, let
$E$ be a set in $\overline{D^{\,\prime}},$ which is closed in
$\overline{{\Bbb R}^n}$ and such that $\partial D^{\,\prime}\subset
E.$ Assume that the following conditions hold:

\medskip
1) for each point $y_0\in \partial D^{\,\prime}\setminus\{\infty\}$
there is $0<r_0:=\sup\limits_{y\in D^{\,\prime}}|y-y_0|$ and any
$0<r_1<r_2<r_0:=\sup\limits_{y\in D^{\,\prime}}|y-y_0|$ there exists
a set $E_1\subset[r_1, r_2]$ of a positive linear Lebesgue measure
is that the function $Q$ is integrable on $S(y_0, r)$ for each $r\in
E_1;$

\medskip
2) for each point $z_0\in E$ and for any neighborhood $U$ of this
point, there is a neighborhood $V\subset U$ of the point $z_0$ such
that the set $V\cap (D^{\,\prime }\setminus E)$ consists of a finite
number of components;

\medskip
3) each component $K$ of the set $D^{\,\prime}\setminus E$ contains
exactly two points $x_K$ and $y_K$ of the set~$P.$

\medskip
Then any $f\in{\frak S}^P_{E, \delta, Q}(D, D^{\,\prime})$ has
continuous extension to mapping
$\overline{f}:\overline{D}\rightarrow \overline{D^{\,\prime}}$ such
that $\overline{f}(\overline{D})=\overline{D^{\,\prime}}$ and the
family ${\frak S}^P_{E, \delta, Q}(D, D^{\,\prime})$ which consists
of all extended mappings $\overline{f}:\overline{D}\rightarrow
\overline{D^{\,\prime}}$ is equicontinuous in $\overline{D}.$ }
\end{theorem}

\section{Preliminaries}

Let $D\subset {\Bbb R}^n,$ $f:D\rightarrow {\Bbb R}^n$ be a discrete
open mapping, $\beta: [a,\,b)\rightarrow {\Bbb R}^n$ be a path, and
$x\in\,f^{\,-1}(\beta(a)).$ A path $\alpha: [a,\,c)\rightarrow D$ is
called a {\it maximal $f$-lifting} of $\beta$ starting at $x,$ if
$(1)\quad \alpha(a)=x\,;$ $(2)\quad f\circ\alpha=\beta|_{[a,\,c)};$
$(3)$\quad for $c<c^{\prime}\leqslant b,$ there is no a path
$\alpha^{\prime}: [a,\,c^{\prime})\rightarrow D$ such that
$\alpha=\alpha^{\prime}|_{[a,\,c)}$ and $f\circ
\alpha^{\,\prime}=\beta|_{[a,\,c^{\prime})}.$ Here and in the
following we say that a path $\beta:[a, b)\rightarrow
\overline{{\Bbb R}^n}$ converges to the set $C\subset
\overline{{\Bbb R}^n}$ as $t\rightarrow b,$ if $h(\beta(t),
C)=\sup\limits_{x\in C}h(\beta(t), C)\rightarrow 0$ at $t\rightarrow
b.$ The following is true (see~\cite[Lemma~3.12]{MRV$_2$}).

\medskip
\begin{proposition}\label{pr3}
{\it\, Let $f:D\rightarrow {\Bbb R}^n,$ $n\geqslant 2,$ be an open
discrete mapping, let $x_0\in D,$ and let $\beta: [a,\,b)\rightarrow
{\Bbb R}^n$ be a path such that $\beta(a)=f(x_0)$ and such that
either $\lim\limits_{t\rightarrow b}\beta(t)$ exists, or
$\beta(t)\rightarrow \partial f(D)$ as $t\rightarrow b.$ Then
$\beta$ has a maximal $f$-lifting $\alpha: [a,\,c)\rightarrow D$
starting at $x_0.$ If $\alpha(t)\rightarrow x_1\in D$ as
$t\rightarrow c,$ then $c=b$ and $f(x_1)=\lim\limits_{t\rightarrow
b}\beta(t).$ Otherwise $\alpha(t)\rightarrow \partial D$ as
$t\rightarrow c.$}
\end{proposition}

\medskip
The following statement plays an important role in the proof of the
main result.

\medskip
\begin{lemma}\label{lem1}
{\,\it Let $D^{\,\prime}\subset {\Bbb R}^n,$ $n\geqslant 2,$ be a
domain in ${\Bbb R}^n,$ let $A\subset \overline{D^{\,\prime}}$ be
such that $A\cap D^{\,\prime}$ is closed with respect to
$D^{\,\prime}.$ Suppose that, $z_0\in A$ and the following is true:
for any neighborhood $U$ of the point $z_0$, there is a neighborhood
$V\subset U$ of the point $z_0$ such that the set $V$ satisfies the
conditions:

\medskip
$1^{*}$) $V\cap (D^{\,\prime}\setminus A)\ne\varnothing,$ and

\medskip
$2^{*}$) $V\cap (D^{\,\prime}\setminus A)$ consists of $m$
components, $m\in {\Bbb N}.$

\medskip
Then, for any neighborhood $U$ of $z_0$ there is an open set
$V\subset U$ such that $z_0\in V$ and

\medskip
1) $V\cap (D^{\,\prime}\setminus A)\ne\varnothing$ and

\medskip
2) $V\cap (D^{\,\prime}\setminus A)$ consists of $l$ components,
where $1\leqslant l\leqslant m.$

\medskip
In other words, <<neighborhood $V$ of the point $z_0$>> may be
replaced by <<open set $V$>>, while the number of corresponding
components does not increase. }
\end{lemma}

\medskip
\begin{proof} Let $U$ be a neighborhood of $z_0.$ Then ${\rm
Int\,}U$ is also neighborhood of $z_0.$ By the assumption, there
exists a neighborhood $V$ of $z_0,$ $V\subset {\rm Int\,}U,$ a
number $m\in {\Bbb N}$ such that
$$V\cap (D^{\,\prime}\setminus A)=V_1\cup V_2\cup\ldots\cup
V_m\ne\varnothing\,.$$
We denote by $W_k,$ $1\leqslant k\leqslant m,$ the corresponding
component of ${\rm Int\,}(U)\cap (D^{\,\prime}\setminus A),$ which
contains the set~$V_k,$ $1\leqslant k\leqslant m$ (some of them may
coincide, but their total number may not exceed $m.$ Let $l$ be the
number of different $W_k$). Since the space $\overline{{\Bbb R}^n}$
is locally path connected and ${\rm Int\,}(U)\cap
(D^{\,\prime}\setminus A)$ is open, then each of the sets $W_k$ is
also open (see~\cite[Theorem~4.II.49.6]{Ku}). Let
$\bigcup\limits_{s=1}^{\infty} W^{\,*}_s$ be the union of all
components of ${\rm Int\,}(U)\cap (D^{\,\prime}\setminus A)$
excluding $\bigcup\limits_{k=1}^mW_k.$ Let us show that the set
\begin{equation}\label{eq7A}
V_*:={\rm Int\,}(U)\setminus \overline{\bigcup\limits_{s=1}^{\infty}
W^{\,*}_s}
\end{equation}
is desired: it is open in $\overline{{\Bbb R}^n},$ it contains $z_0$
and $V_*\cap (D^{\,\prime}\setminus A)$ consists of exactly $l$
component.

The openness of $V_*$ is obvious, since $V_*$ is the difference of
two sets, the first of which is open, and the second is closed.
Next, note that the set $\bigcup\limits_{s=1}^{\infty} W^{\,*}_s$ is
closed in ${\rm Int\,}(U)\cap (D^{ \,\prime}\setminus A).$ Then
\begin{eqnarray*}V_*\cap (D^{\,\prime}\setminus A)={\rm
Int\,}(U)\cap\left(D^{\,\prime}\setminus A\right) \setminus
\overline{\bigcup\limits_{s=1}^{\infty} W^{\,*}_s}\\={\rm
Int\,}(U)\cap\left(D^{\,\prime}\setminus A\right) \setminus
\bigcup\limits_{s=1}^{\infty} W^{\,*}_s=W_1\cup W_2\cup\ldots\cup
W_l\,. \end{eqnarray*}
Finally, it is not difficult to prove from the opposite that $z_0\in
V_*.$ Indeed, if $z_0\not\in V_*,$ then there is a sequence $z_k\in
\bigcup\limits_{s=1}^{\infty} W^{\,*}_s$ such that $z_k\rightarrow
z_0$ as $k\rightarrow\infty.$ In particular, $z_k\in
D^{\,\prime}\setminus A.$ On the other hand, since $V$ is a
neighborhood of $z_0,$ for large $k\in {\Bbb N}$ we have that
$$z_k\in V\cap (D^{\,\prime}\setminus A)=V_1\cup V_2\cup\ldots\cup
V_m\subset\bigcup\limits_{k=1}^mW_k,$$
but this contradicts the assumption that
$z_k\in\bigcup\limits_{s=1}^{\infty} W^{\,*}_s,$ because by the
definition of components,
$$\left(\bigcup\limits_{s=1}^{\infty} W^{\,*}_s\right)\cap\left(
\bigcup\limits_{k=1}^mW_k\right)=\varnothing\,.$$ The lemma is
proved.~$\Box$
\end{proof}

\section{Proof of Theorem~\ref{th3}}

We fix $x_0\in\partial D.$ Let us show that the mapping $f$ has a
continuous extension to a point $x_0.$ Using a M\"{o}bius
transformation $\varphi:\infty\mapsto 0$ if necessary, due to the
invariance of the modulus $M$ in~(\ref{eq2*A})
(see~\cite[Theorem~8.1]{Va$_2$}), we may assume that $x_0\ne\infty.$

\begin{figure}[h]
\centerline{\includegraphics[scale=0.5]{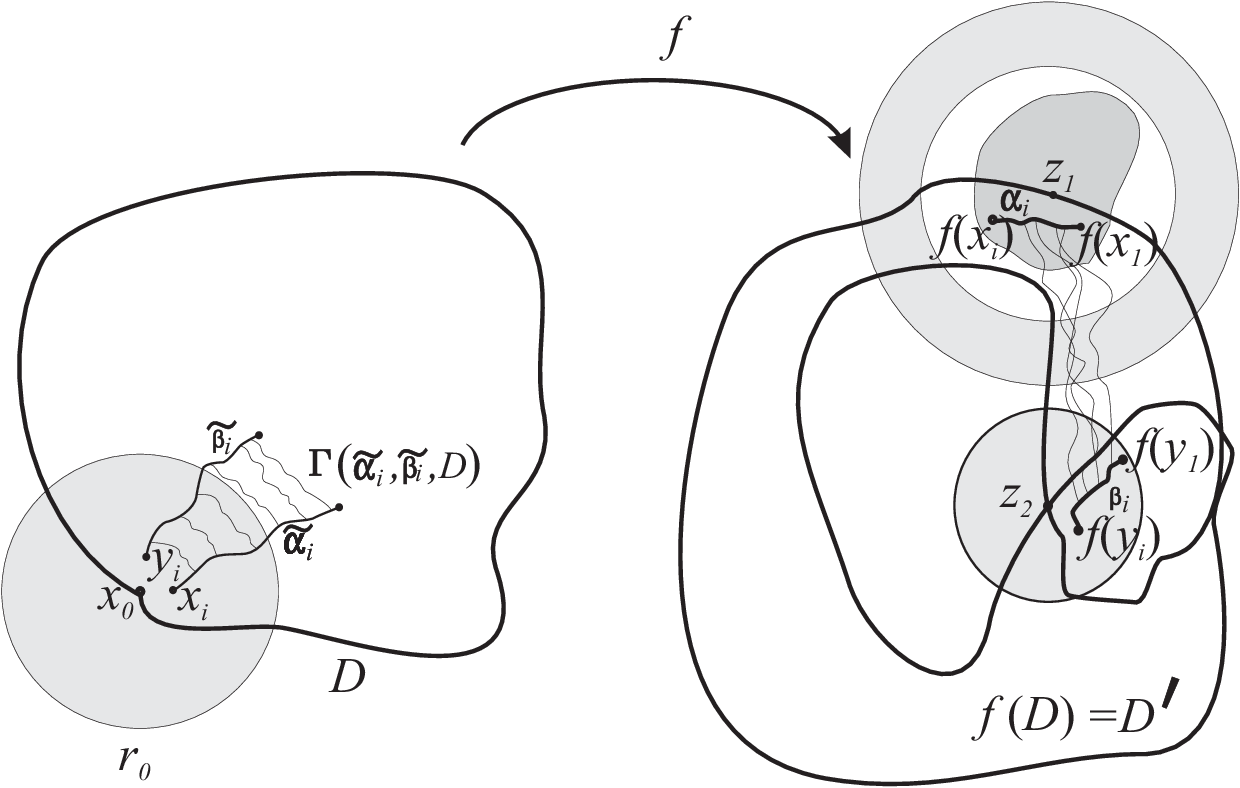}} \caption{To
the proof of Theorem~\ref{th3}}\label{fig1}
\end{figure}

\medskip
Let us carry out the proof from the opposite, namely, suppose that
$f$ has no a continuous extension to the point $x_0.$ Then there
exist sequences $x_i, y_i\in D,$ $i=1,2,\ldots ,$ such that $x_i,
y_i\rightarrow x_0$ as $i\rightarrow\infty,$ and
\begin{equation}\label{eq1B}
h(f(x_i), f(y_i))\geqslant a>0
\end{equation}
for some $a>0$ and all $i\in {\Bbb N},$ where $h$ is the chordal
(spherical) metric, see~(\ref{eq3C}). Due to the compactness of
$\overline{{\Bbb R}^n},$ we may assume that the sequences $f(x_i)$
and $f(y_i)$ converge as $i\rightarrow\infty$ to $z_1$ and $z_2,$
respectively, and $z_1\ne\infty,$ see Figure~\ref{fig1}.
First of all, we note that the points $x_i$ and $y_i,$
$i=1,2,\ldots, $ may be chosen such that $x_i, y_i\not\in
f^{\,-1}(E\cap D^{\,\prime}).$ Indeed, since under condition~3) the
set $f^{\,-1}(E\cap D^{\,\prime})$ is nowhere dense in $D,$ there
exists a sequence $x_{ki}\in D\setminus(f^{\,-1}(E\cap
D^{\,\prime})),$ $k=1,2,\ldots ,$ such that $x_{ki}\rightarrow x_i$
as $k\rightarrow\infty.$ Put $\varepsilon>0.$ Due to the continuity
of the mapping $f$ at the point $x_i,$ for the number $i\in {\Bbb
N}$ there is a number $k_i\in {\Bbb N}$ such that $h(f(x_{k_ii}),
f(x_i))<\frac{1}{2^i}.$ So, by the triangle inequality
$$h(f(x_{k_ii}), z_1)\leqslant h(f(x_{k_ii}), f(x_i))+
h(f(x_i), z_1)\leqslant \frac{1}{2^i}+\varepsilon\,,$$
$i\geqslant i_0=i_0(\varepsilon),$ since $f(x_i)\rightarrow z_1$ as
$i\rightarrow\infty$ and by the choice of $x_i$ and $z_1.$
Therefore, $x_i\in D$ may be replaced by $x_{k_ii}\in
D\setminus(f^{\,-1}(E\cap D^{\,\prime})),$ as required. We may
reason similarly for the sequence $y_i.$

\medskip
Since $E$ is a closed set in $\overline{{\Bbb R}^n},$ then, in
particular, $E$ is closed with respect to $D^{\,\prime}.$ Then the
set $D^{\,\prime}\setminus E=D^{\,\prime}\setminus (E\cap
D^{\,\prime})$ is open. Next, by condition~2), given neighborhoods
$U_1$ and $U_2$ of the points $z_1$ and $z_2\in E,$ there are
neighborhoods $V_1\subset U_1$ and $V_2\subset U_2$ of these same
points such that each of the sets $W_1:=V_1\cap
(D^{\,\prime}\setminus E)$ and $W_2:=V_2\cap (D^{\,\prime}\setminus
E)$ consist of a finite number of components. By Lemma~\ref{lem1},
$V_1$ and $V_2$ may be chosen to be open. Then $W_1$ and $W_2$ are
also open, as well. Therefore, each component of $W_1$ and $W_2$ is
open (see~\cite[Theorem~4.II.49.6]{Ku}) and path connected (see,
e.g.,~\cite[Proposition~13.2]{MRSY$_2$}).

\medskip
Without loss of generality, we may assume that the number
$r_0=r_0(z_1)>0$ and the numbers $0<r_1<r_2<r_0$ from the condition
of the theorem are such that
\begin{equation}\label{eq8}
U_1\subset B(z_1, r_1),\qquad \overline{B(z_1, r_2)}\cap
\overline{U_2}=\varnothing\,.
\end{equation}
Passing to subsequences if necessary and taking into account that
$W_1$ and $W_2$ have a finite number of components, we may assume
that $f(x_i),$ $i=1,2,\ldots ,$ belong to some single component
$K_1$ of the set $W_1.$ Similarly, we may assume that $f(y_i),$
$i=1,2,\ldots ,$ belong to single component $K_2$ of the set $W_2.$
Let us to join the points $f(x_i)$ and $f(x_1)$ by a path
$\alpha_i:[0, 1]\rightarrow D^{\,\prime},$ and the points $f(y_i)$
and $f(y_1)$ by a path $\beta_i:[0, 1]\rightarrow D^{\,\prime}$ so
that $|\alpha_i|\subset K_1$ and $|\beta_i|\subset K_2$ for
$i=1,2,\ldots .$ Let $\widetilde{\alpha_i}:[0, c_1)\rightarrow D$
and $\widetilde{\beta_i}:[0, c_2)\rightarrow D^{\,\prime}$ be
maximal $f$-liftings of paths $\alpha_i$ and $\beta_i$ starting at
$x_i$ and $y_i,$ respectively (they exist by Proposition~\ref{pr3}).
Due to the same proposition, only one of the following two
situations are possible:

\medskip
1) $\widetilde{\alpha_i}(t)\rightarrow x_1\in D$ as $t\rightarrow
c_1,$ and $c_1=1$ and $f(\widetilde{\alpha_i}(1))=f(x_1),$ or

\medskip
2) $\widetilde{\alpha_i}(t)\rightarrow \partial D$ as $t\rightarrow
c_1.$

\medskip
Let us to show that situation 2) is impossible. Suppose the
opposite: let $\widetilde{\alpha_i}(t)\rightarrow \partial D$ as
$t\rightarrow c_1.$ We choose an arbitrary sequence $t_m\in [0,
c_1)$ such that $t_m\rightarrow c_1-0$ as $m\rightarrow\infty.$
Since the space $\overline{{\Bbb R}^n}$ is compact, the boundary
$\partial D$ is also compact as a closed subset of the compact
space. Then there exists $w_m\in \partial D$ such that
\begin{equation}\label{eq7B}
h(\widetilde{\alpha_i}(t_m), \partial
D)=h(\widetilde{\alpha_i}(t_m), w_m) \rightarrow 0\,,\qquad
m\rightarrow \infty\,.
\end{equation}
Due to the compactness of $\partial D$, we may assume that
$w_m\rightarrow w_0\in \partial D$ as $m\rightarrow\infty.$
Therefore, by the relation~(\ref{eq7B}) and by the triangle
inequality
\begin{equation}\label{eq8B}
h(\widetilde{\alpha_i}(t_m), w_0)\leqslant
h(\widetilde{\alpha_i}(t_m), w_m)+h(w_m, w_0)\rightarrow 0\,,\qquad
m\rightarrow \infty\,.
\end{equation}
On the other hand,
\begin{equation}\label{eq9B}
f(\widetilde{\alpha_i}(t_m))=\alpha_i(t_m)\rightarrow
\alpha_i(c_1)\in D^{\,\prime}\setminus E\,,\quad
m\rightarrow\infty\,,
\end{equation}
because by the construction the path $\alpha_i(t),$ $t\in [0, 1],$
lies in $D^{\,\prime}\setminus E$ together with its finite ones
points. At the same time, by~(\ref{eq8B}) and~(\ref{eq9B}) we have
that $\alpha_i(c_1)\in C(f,
\partial D)\subset E.$ The inclusion $|\alpha_i|\subset D^{\,\prime}\setminus
E$ contradicts the relation $\alpha_i(c_1)\in E,$ where, as usual,
given a path $\gamma:[a, b]\rightarrow \overline{{\Bbb R}^n}$ we use
the notation
$$|\gamma|=\{x\in\overline{{\Bbb R}^n}: \exists\,t\in [a,
b]:\gamma(t)=x\}\,$$ for the {\it locus} of $\gamma.$ The resulting
contradiction indicates the impossibility of the second situation
mentioned above.

\medskip
Therefore, the first situation is fulfilled:
$\widetilde{\alpha_i}(t)\rightarrow x_1\in D$ as $t\rightarrow c_1,$
and $c_1=1$ and $f(\widetilde{\alpha_i}(1))=f(x_1).$ In other words,
the $f$-lifting $\widetilde{\alpha_i}$ is complete, i.e.,
$\widetilde{\alpha_i}:[0, 1]\rightarrow D.$ Similarly, the path
$\beta_i$ has a complete $f$-lifting $\widetilde{\beta_i}:[0,
1]\rightarrow D.$

\medskip
Let us to prove that, the points $f(x_1)$ and $f(y_1)$ have only a
finite number of pre-images in $D.$ Let us prove this from the
opposite: suppose that there is a sequence $z_i,$ $i=1,2,\ldots ,$
such that $f(z_i)=f(x_1).$ Due to the compactness of
$\overline{{\Bbb R}^n}$ we may assume that $z_i$ converges to some
point $z_0$ as $i\rightarrow\infty.$ Since $f$ is discrete, $z_0\in
\partial D.$ But then $f(x_1)\in C(f, \partial D)\subset E,$ which contradicts
the definition of $f(x_1).$ The resulting contradiction indicates on
the finiteness of number of pre-images of the point $f(x_1)$ in $D$
under the mapping $f.$ The similar conclusion may be done with
respect to $f(y_1).$

\medskip
Therefore, there exists $\widetilde{r}_0>0$ such that
$\widetilde{\alpha_i}(1), \widetilde{\beta_i}(1)\in D\setminus
B(x_0, \widetilde{r}_0)$ for all $i=1,2,\ldots .$ Since the boundary
of the domain $D$ is weakly flat, for every $P>0$ there is
$i=i_P\geqslant 1$ such that
\begin{equation}\label{eq7}
M(\Gamma(|\widetilde{\alpha_i}|, |\widetilde{\beta_i}|,
D))>P\qquad\forall\,\,i\geqslant i_P\,.
\end{equation}
Let us to show that the condition~(\ref{eq7}) contradicts the
definition of the mapping $f$ in~(\ref{eq2*A}). Indeed,
by~(\ref{eq8}) and by~\cite[Theorem~1.I.5.46]{Ku}
\begin{equation}\label{eq9}
f(\Gamma(|\widetilde{\alpha_i}|, |\widetilde{\beta_i}|,
D))>\Gamma(S(z_1, r_1), S(z_1, r_2), A(z_1, r_1, r_2))\,.
\end{equation}
Therefore, it follows from~(\ref{eq9}) that
\begin{equation}\label{eq10}
\Gamma(|\widetilde{\alpha_i}|, |\widetilde{\beta_i}|, D)
>\Gamma_f(z_1, r_1, r_2)\,.
\end{equation}
In turn, by~(\ref{eq10}) we have that
\begin{equation}\label{eq11}
M(\Gamma(|\widetilde{\alpha_i}|, |\widetilde{\beta_i}|, D))\leqslant
M(\Gamma_f(z_1, r_1, r_2))\leqslant \int\limits_{A} Q(y)\cdot \eta^n
(|y-z_1|)\, dm(y)\,,
\end{equation}
where $A=A(z_1, r_1, r_2)$ and $\eta$ is an arbitrary Lebesgue
measurable function satisfying the relation~(\ref{eqA2}). We use the
following standard conventions: $a/\infty=0$ for $a\ne\infty,$
$a/0=\infty$ for $a>0$ and $0\cdot\infty=0$ (see, e.g.,
\cite[3.I]{Sa}). Set $\widetilde{Q}(y)=\max\{Q(y), 1\},$ and let
\begin{equation}\label{eq13}
I=\int\limits_{r_1}^{r_2}\frac{dt}{t\widetilde{q}_{z_1}^{1/(n-1)}(t)}\,,
\end{equation}
where
\begin{equation}\label{eq12}
\widetilde{q}_{z_1}(r)=\frac{1}{\omega_{n-1}r^{n-1}}\int\limits_{S(z_1,
r)}\widetilde{Q}(y)\,d\mathcal{H}^{n-1}(y)\,, \end{equation}
and $\omega_{n-1}$ is the $(n-1)$-dimensional area of the unit
sphere ${\Bbb S}^{n-1}$ in ${\Bbb R}^n.$ By assumption, there is a
set $E_1\subset [r_1 , r_2]$ of positive one-dimensional Hausdorff
measure such that $q_{z_1}(t)$ is finite at of all $t\in E_1.$
Consequently, $\widetilde{q}_{z_1}(t)$ is finite at of all $t\in
E_1,$ as well. Therefore, $I\ne 0$ in~(\ref{eq13}). Besides that,
$I\ne\infty,$ because
$$I\leqslant \log\frac{r_2}{r_1}<\infty\,.$$
Put $\eta_0(t)=(Itq_{z_1}^{1/(n-1)}(t))^{\,-1}$ and observe that
$\eta_0$ satisfies the relation~(\ref{eqA2}). Let us substitute this
function in the right part of the inequality~(\ref{eq11}) and apply
Fubini's theorem. We will have that
\begin{equation}\label{eq14}
M(\Gamma(|\widetilde{\alpha_i}|, |\widetilde{\beta_i}|, D))
\leqslant \frac{\omega_{n-1}}{I^{n-1}}<\infty\,.
\end{equation}
The relation~(\ref{eq14}) contradicts~(\ref{eq7}). The contradiction
obtained above shows that the assumption in~(\ref{eq1B}) is not
true.

\medskip
Let $\overline{f}$ be a continuous extension of $f$ onto
$\overline{D}.$ Now, we show that
$\overline{f}(D)=\overline{D^{\,\prime}}.$ Indeed, by the definition
of mapping $\overline{f}$ we immediately have that
$\overline{f}(\overline{D})\subset\overline{D^{\,\prime}}.$ Let us
show that $\overline{D^{\,\prime}}\subset
\overline{f}(\overline{D}).$ Let $y_0\in \overline{D^{\,\prime}},$
then either $y_0\in D^{\,\prime},$ or $y_0\in \partial
D^{\,\prime}.$ If $y_0\in D^{\,\prime},$ then
$y_0=f(\widetilde{x}_0)$ and $y_0\in \overline{f}(\overline{D}),$
because $f$ maps $D$ onto $D^{\,\prime}$ by the assumption. Finally,
let $y_0\in
\partial D^{\,\prime}.$ Then there is a sequence $y_i\in D^{\,\prime}$
such that $y_i=f(x_i)\rightarrow y_0$ as $i\rightarrow\infty$ and
$x_i\in D.$ Since the space $\overline{{\Bbb R}^n}$ is compact, we
may assume that $x_i\rightarrow \widetilde{x}_0,$ where
$\widetilde{x}_0\in\overline{D}.$ Note that, $\widetilde{x}_0\in
\partial D,$ since $f$ is open. Then $f(\widetilde{x}_0)=y_0\in
\overline{f}(\partial D)\subset \overline{f}(\overline{D}).$
Theorem~\ref{th3}.~$\Box$

\medskip
{\it Proof of Corollary~\ref{cor1}. Let $r_0:=\sup\limits_{y\in
D^{\,\prime}}|y-y_0|.$ We may assume that $Q$ is extended by zero
outside $D^{\,\prime}.$ By the Fubini theorem (see, e.g.,
\cite[Theorem~8.1.III]{Sa}) we obtain that
$$\int\limits_{r_1<|y-y_0|<r_2}Q(y)\,dm(y)=\int\limits_{r_1}^{r_2}
\int\limits_{S(y_0, r)}Q(y)\,d\mathcal{H}^{n-1}(y)dr<\infty\,.$$
This means the fulfillment of condition~1) in
Theorem~\ref{th3}.~$\Box$}

\section{Some examples}

\begin{example}\label{ex4}
Consider the disk $D=B(1, 1)=\{z\in {\Bbb C}: |z-1|<1\}$ on the
complex plane. Let us define a mapping $f:D\rightarrow {\Bbb C}$ as
follows: $f(z)=z^4.$ The Figure~\ref{fig2} schematically shows the
image of the domain $D$ under the mapping $f;$ the complete image is
marked in gray, the domain that is mapped twice in the disk $B(1,
1)$ is marked in a darker color.
\begin{figure}[h]
\centerline{\includegraphics[scale=0.4]{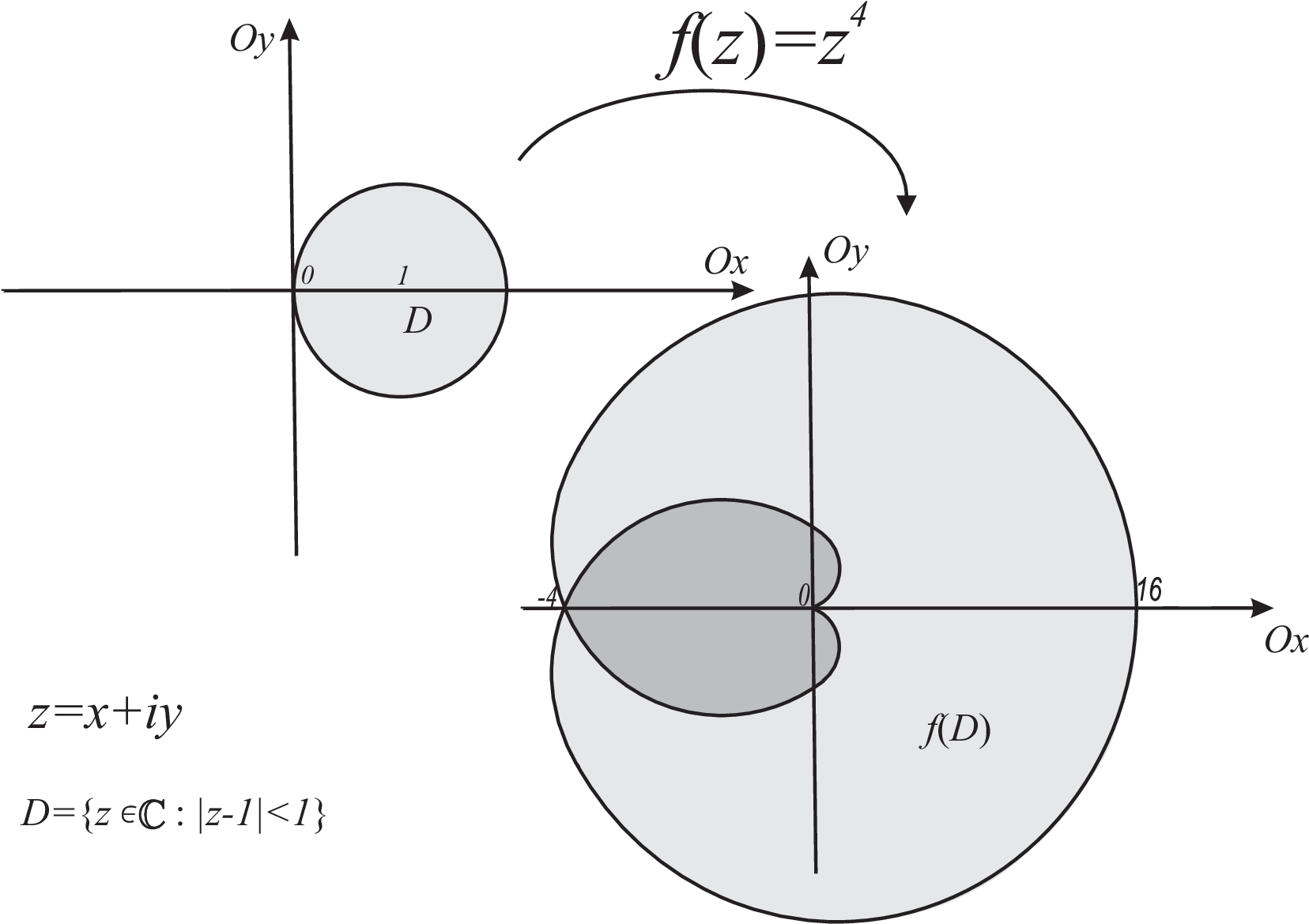}}
\caption{A non-boundary preserving quasiregular mapping that
satisfies conditions of Theorem~\ref{th3} }\label{fig2}
\end{figure}
The image of the boundary of the domain $D$ has one point of
self-intersection $w=(-4, 0).$ Paths
$z_1(\varphi)=16\cos^4\varphi\cdot e^{i\cdot4\varphi},$
$\frac{\pi}{4}<\varphi \leqslant \frac{\pi}{2},$ and
$z_2(\varphi)=16\cos^4\varphi\cdot e^{i\cdot 4\varphi},$
$-\frac{\pi}{2}\leqslant\varphi<-\frac{\pi}{4},$ is a part of the
image of the boundary of $D$ that lies inside the mapped domain. It
is not difficult to see that the complete pre-image in $D$ of this
paths consists of two paths
$\widetilde{z_1}(\varphi)=2\cos\varphi\cdot
e^{i\left(\varphi+\frac{\pi}{2}\right)},$ $\frac{\pi}{4}<\varphi
\leqslant \frac{\pi}{2},$ and
$\widetilde{z_2}(\varphi)=2\cos\varphi\cdot
e^{i\left(\varphi-\frac{\pi}{2}\right)},$
$-\frac{\pi}{2}\leqslant\varphi<-\frac{\pi}{4}.$

Note that, the mapping $f$ satisfies all conditions of
Theorem~\ref{th3} with $E=C(f, \partial D).$ Indeed, $f(D)$ is
locally finitely connected with respect to $C(f, \partial D)$
(simply connected at points of outer boundary, excluding the point
of self-intersection $w=(-4, 0);$ is three-connected at this point
and two-connected at all other points of $C(f, \partial D)$). In
addition, due to the mentioned above, $f^{\,-1}(C(f, \partial D)\cap
f(D))$ is nowhere dense in $D,$ since $f^{\,-1}(C(f, \partial D)\cap
f(D))$ consists of two paths $\widetilde{z_1}(\varphi)$ and
$\widetilde{z_2}(\varphi),$ whose closure is nowhere dense in $D.$
Finally, $f(z)=z^4$ is a quasiregular mapping, therefore satisfies
the relation~(\ref{eq2*A}) for the function $Q(y)\equiv 4$ (see,
e.g., \cite[Theorem~3.2]{MRV$_1$} and Remark~\ref{rem1}). Obviously,
the function $Q(z)=4$ is integrable over almost all circles $S(z_0,
r)$. (Moreover, it is integrable in $D$ with respect to Lebesgue
measure). By Theorem~\ref{th3} $f$ has a continuous extension to
$\partial D.$
\end{example}

\medskip
\begin{example}\label{ex5}
It is not difficult to construct a similar mapping with some
unbounded function $Q$ in~(\ref{eq2*A}). For this purpose, consider
an arbitrary point $z_0\in f(D)\setminus f(\partial D)$ in the
previous example, and any $0<r_0<{\rm dist\,}(z_0, f(\partial D)).$
Let
$$h(z)=\frac{r_0z}{|z|\log\frac{e}{|z|}}+z_0\,,\qquad h(0)=z_0\,,$$
$g(y):=h^{\,-1}(y).$ Then $g$ is well-defined in the disk $B(z_0,
r_0)$ and $g(B(z_0, r_0))={\Bbb D},$ ${\Bbb D}=\{z\in {\Bbb C}:
|z|<1\}.$ Reasoning similarly to~\cite[Proposition~6.3]{MRSY$_2$},
we may show that $g$ satisfies the
relation~(\ref{eq2*A})--(\ref{eqA2}) in each point $y_0\in
\overline{{\Bbb B}^n}$ for $Q=Q(y)=\log\left(\frac{e}{|y|}\right).$

Note that, $Q\in L^1({\Bbb D}).$ Indeed, by the Fubini theorem, we
will have
\begin{eqnarray*}\int\limits_{{\Bbb D}}Q(y)\,dm(y)=\int\limits_0^1\int\limits_{S(0,
r)} \log\left(\frac{e}{|y|}\right)d\mathcal{H}^{1}(y)dr \\
=2\pi\int\limits_0^1 r\log\left(\frac{e}{r}\right)\,dr\leqslant 2\pi
e\int\limits_0^1dr=2\pi e<\infty\,. \end{eqnarray*}
Therefore, the mapping $g:B(z_0, r_0)\rightarrow {\Bbb C}$ satisfies
all conditions of Theorem~\ref{th3}.

Note that, $h(z)=r_0z+z_0$ for $z\in {\Bbb S},$ ${\Bbb S}=\{z\in
{\Bbb C}: |z|=1\}.$ Then $z=(h-z_0)/r_0,$ $g(y)=(y-z_0)/r_0.$ Put
$$G(z)=\begin{cases}(g\circ f)(z)\,,& z\in f^{\,-1}(B(z_0, r_0))\,, \\
(f(z)-z_0)/r_0\,, & z\not\in f^{\,-1}(B(z_0, r_0))
\end{cases}\,.$$
The mapping $G$ is constructed so that it does not change the
geometry of the domain $f(D)$ with accuracy up to the shift to point
$z_0$ for compression $r_0$ times. In addition, the mapping $G$
satisfies the inequality~(\ref{eq2*A})--(\ref{eqA2}) for
$$Q_1(y):=\begin{cases}\log\left(\frac{e}{|y|}\right)\,,& y\in {\Bbb
D}\,,\\ 1\,,& y\in {\Bbb C}\setminus {\Bbb D}\,,
\end{cases}$$
and, moreover,
\begin{equation}\label{eq25}
f(\Gamma_{G}(y_0, r_1, r_2))\subset \Gamma_{g}(y_0, r_1, r_2)\,.
\end{equation}
Obviously, $N(f, D)\leqslant 4.$ Then, due to~(\ref{eq25}) and
Remarks~\ref{rem1} we obtain that
\begin{multline}\label{eq26}
M(\Gamma_{G}(y_0, r_1, r_2))\leqslant 4\cdot M(f(\Gamma_{g}(y_0,
r_1, r_2)))\leqslant 4\cdot M(\Gamma_{g}(y_0, r_1, r_2)) \\
\leqslant 4\cdot \int\limits_{(G(D)\cap A(y_0,r_1,r_2)} Q_1(y)\cdot
\eta^2 (|y-y_0|)\, dm(y)\,,
\end{multline}
where $\eta: (r_1,r_2)\rightarrow [0,\infty ]$ is any Lebesgue
measurable function such that
$$
\int\limits_{r_1}^{r_2}\eta(r)\, dr\geqslant 1\,.
$$
A mapping $G$ is open and discrete as a composition of a
quasiregular mapping $f$ with some homeomorphism $G.$ Now, by
Theorem~\ref{th3} it has a continuous extension in $\overline{D}.$
\end{example}

\section{On equicontinuity in the closure of a domain}

A statement similar in content to the one below has been published
in~\cite[Lemma~2.1]{SevSkv$_2$}.

\medskip
\begin{lemma}\label{lem3}
{\it\, Let $D$ be a domain in ${\Bbb R}^n,$ $n\geqslant 2,$ and let
the points $a\in D, b\in \overline{D},$ and $c\in D, d\in
\overline{D}$ be such that there are paths $\gamma_1:[0,
1]\rightarrow \overline{D}$ and $\gamma_2:[0, 1]\rightarrow
\overline{D},$ such that $\gamma_i(t)\in D$ for all $t\in (0, 1),$
$i=1,2,$ $\gamma_1(0)=a,$ $\gamma_1(1)=b,$ $\gamma_2(0)= c,$
$\gamma_2(1)=d.$ Then there are disjoint paths with the same
properties $\gamma^{\,*}_1$ and $\gamma^{\,*}_2.$

Moreover, if $\gamma_1$ is Jordanian, then we may take $\gamma_1$ as
$\gamma^{\,*}_1,$ and as $\gamma^{\,*}_2$ we may take a path which
coincides with $\gamma_2$ on some segments $[0, t_0]$ and $[T_0,
1],$ where $0<t_0<T_0<1.$}
\end{lemma}

\medskip
\begin{proof}
Discarding, if necessary, a (not more than countable) number of
loops in $\gamma_1,$ we may consider that $\gamma_1$ is Jordanian.
Set now $\gamma^{\,*}_1:=\gamma_1.$

Let $n\geqslant 3.$ Then $\gamma_1$ does not split $D$ as a set of
topological dimension~1 (see \cite[Consequence~1.5.IV]{HW}). Since
$\gamma_2(0)=c$ and $\gamma_2(1)=d,$ there exist $0<t_0<T_0<1$ such
that $\gamma_2(t)\not\in |\gamma_1|$ for all $0\leqslant t\leqslant
t_0$ and $T_0\leqslant t\leqslant 1.$ Join the points
$\gamma_2(t_0)$ and $\gamma_2(T_0)$ by a path
$\widetilde{\gamma}:[t_0, T_0]\rightarrow D,$ which does not
intersect $\gamma_1$ (this is possible since $\gamma_1$ does not
split $D$). Put
$$\gamma^{\,*}_2(t):=\begin{cases}\widetilde{\gamma}(t)\,,& t\in [t_0, T_0]\,,\\
\gamma_2(t)\,,& t\in [0, t_0]\cup[T_0, 1]\,.\end{cases}$$
The path $\gamma^{\,*}_2$ is desired: it belongs to $D,$ except for
the second endpoint, join the points $c$ and $d$ and does not
intersect with $\gamma_1.$

\medskip
Now let $n=2.$ Due to the Antoine's theorem (see \cite[Theorem~4.3,
\S\,4]{Keld}) the domain $D$ may be mapped onto some domain
$D^{\,*}$ using a plane homeomorphism $\varphi:{\Bbb R}^2\rightarrow
{\Bbb R}^2$ so that $\varphi(\gamma_1)=J$ and $J$ is a segment in
$D^{\,*}.$ Accordingly, the points $\varphi(c)$ and $\varphi(d)$ may
be joined in $D^{\,*}$ by some path $\alpha_2:[0, 1]\rightarrow
\overline{D^ {\,*}},$ which completely lies in $D^{\,*},$ except,
perhaps, for its endpoint $\alpha_2(1)=\varphi(d).$

It remains to show that the path $\alpha_2$ may be chosen so that it
does not intersect the segment $J.$ In fact, let $\alpha_2$
intersect $J,$ and let $t_1$ and $t_2$ be, respectively, the
smallest and the largest values $t\in [0, 1],$ for which
$\alpha_2(t)\in |J|.$ Also let
$$J=J(s)=\varphi(a)+ (\varphi(b)-\varphi(a))s, \quad s\in [0, 1]$$
be a parameterization of the segment $J.$ Let $\widetilde{s_1}$ and
$\widetilde{s_2}\in (0, 1)$ be such that
$J(\widetilde{s_1})=\alpha_2(t_1)$ and
$J(\widetilde{s_2})=\alpha_2(t_2). $ Let $s_2=\max\{\widetilde{s_1},
\widetilde{s_2}\}.$ Let $e_1=\varphi(b)-\varphi(a)$ and let $e_2$ be
a unit vector, orthogonal $e_1.$ Set
$$P_{\varepsilon}=\{x=\varphi(a)+ x_1e_1+x_2e_2,\quad x_1\in (-\varepsilon, s_2+\varepsilon),
\quad x_2\in (-\varepsilon, \varepsilon)\}\,,\quad
\varepsilon>0\,,$$
is a rectangle containing $|J_1|,$ where $J_1$ is a restriction of
$J$ onto the segment $[0, s_2]$ (see Figure~\ref{fig7}).
\begin{figure}[h] \centerline{\includegraphics[scale=0.5]{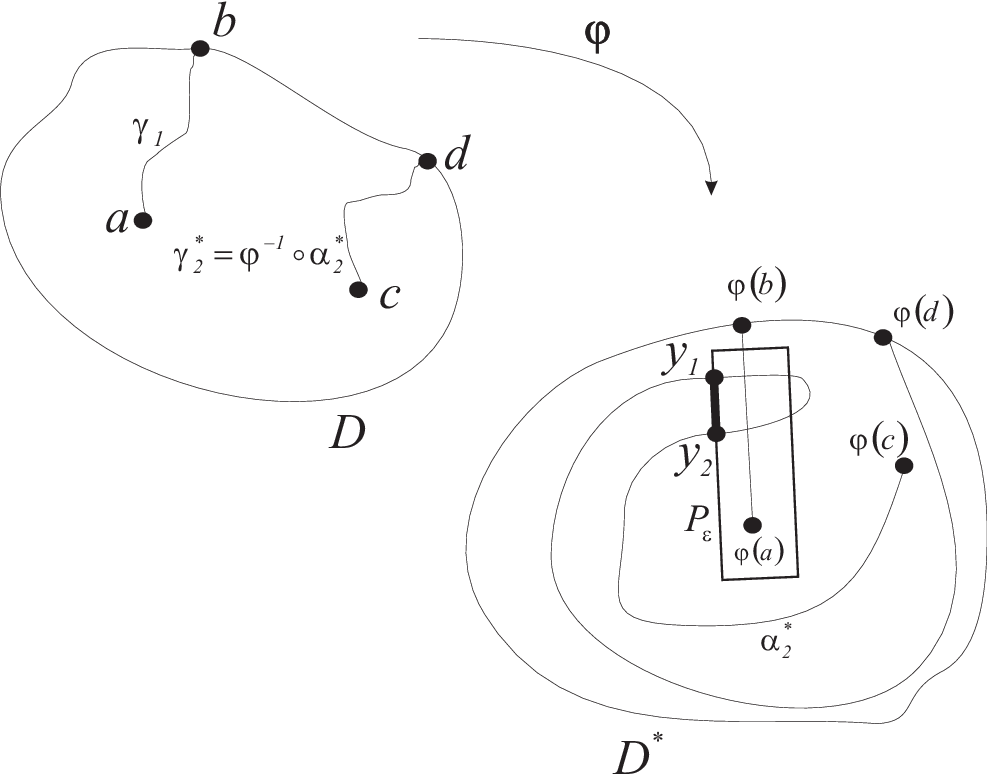}}
\caption{To the proof of Lemma~\ref{lem3}}\label{fig7}
\end{figure}
We choose $\varepsilon>0$ so that $\varphi(c)\not\in
P_{\varepsilon},$ ${\rm dist}\,(P_{\varepsilon}, \partial
D^{\,*})>\varepsilon.$ Due to \cite[Theorem~1.I, ch.~5, \S\,
46]{Ku}) a path $\alpha_2$ intersects $\partial P_{\varepsilon}$ for
some $T_1<t_1$ and $T_2>t_2.$ Let $\alpha_2(T_1)=y_1$ and
$\alpha_2(T_2)=y_2.$ Since $\partial P_{\varepsilon}\setminus
\{z_0\},$ $z_0:=\varphi(a)+(s_2+\varepsilon)e_1,$ is a connected
set, it is possible join the points $y_1$ and $y_2$ with a path
$\alpha^{\,*}(t):[T_1, T_2]\rightarrow \partial
P_{\varepsilon}\setminus \{z_0\}.$ Finally, let
$$\alpha_2^{\,*}(t)\quad =\quad\left\{
\begin{array}{rr} \alpha_2(t), & t\in [0, 1]
\setminus [T_1, T_2],\\ \alpha^{\,*}(t), & t\in [T_1,
T_2]\end{array} \right.$$
and $\gamma^{\,*}_2:=\varphi^{\,-1}\circ \alpha_2^{\,*}.$ Then
$\gamma_1$ joins $a$ and $b$ in $\overline{D},$ and $\gamma_2^{\,*}$
joins $c$ and $d$ in $\overline{D},$ while $\gamma_1$ and $
\gamma_2^{\,*}$ are disjoint, which should have been established.
This time $\gamma^{\,*}_2$ coincides with $\gamma_2$ on the segments
$[0, t_0]$ and $[T_0, 1],$ where $t_0:=T_1$ and $T_0:= T_2.$$~\Box$
\end{proof}

\medskip
The following lemma is an analogue of Lemma~3.10 in~\cite{Vu} on
accessibility of boundary points of domains finitely connected at
the boundary.

\medskip
\begin{lemma}\label{lem4}
{\it\, Let $D^{\,\prime}$ be a domain in ${\Bbb R}^n,$ $n\geqslant
2,$ and let $E\subset \overline{D^{\,\prime}}$ be a closed set in
$\overline{{\Bbb R}^n}.$ Suppose that, for any $z_0\in E$ and for
any neighborhood $U$ of $z_0$ there is a neighborhood $V\subset U$
such that $(V\setminus E)\cap D^{\,\prime}$ consists of a finite
number of components. Let $z_1\in E,$ and let $x_k\in
D^{\,\prime}\setminus E,$ $k=1,2,\ldots ,$ be a sequence that
converges to $z_1$ as $k\rightarrow\infty.$ Then there is s
subsequence $x_{k_l},$ $l=1,2,\ldots ,$ which belongs to exactly one
component $K_1$ of $D^{\,\prime}\setminus E.$

Moreover, if $a_1\in K_1,$ then there is (possibly) a new
subsequence $x_{k_{l_m}},$ $m=1,2,\ldots ,$ and a path $\alpha:[0,
1]\rightarrow \overline{D^{\,\prime}}$ such that:

\medskip
a) $\alpha(t)\in D^{\,\prime}\setminus E$ at $0\leqslant t<1,$

\medskip
b) $x_{k_{l_m}}\in |\alpha|,$ $m=1,2,\ldots ,$

\medskip
c) $\alpha(0)=a_1,$ $\alpha(1)=z_1,$ see Figure~\ref{fig4}.
\begin{figure}[h]
\centerline{\includegraphics[scale=0.4]{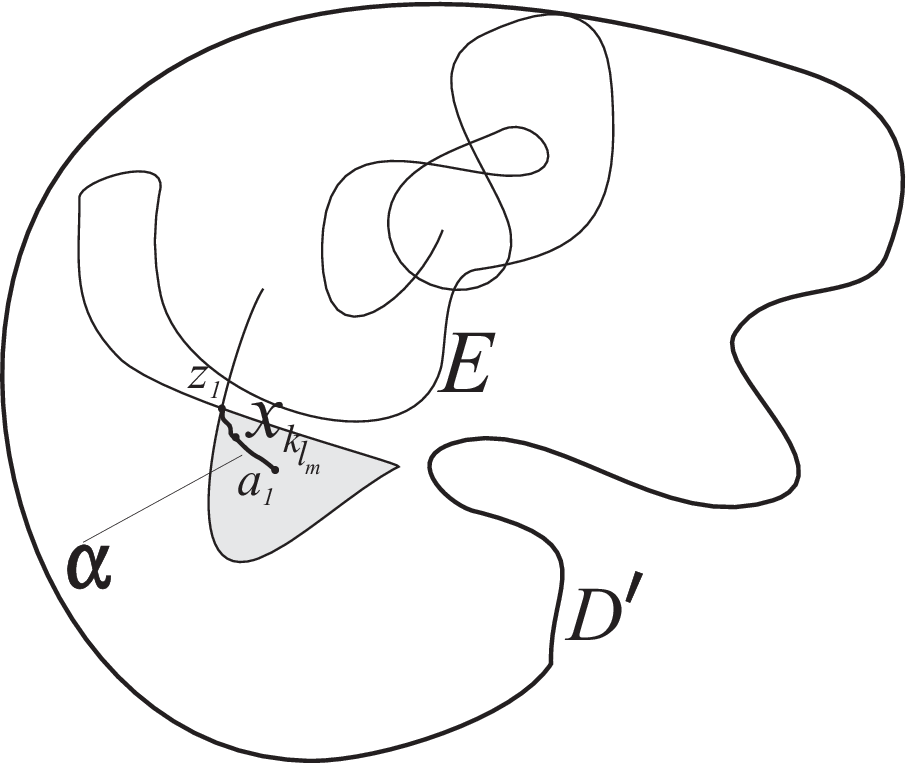}}
\caption{Formulation of Lemma~\ref{lem4}}\label{fig4}
\end{figure}}
\end{lemma}

\medskip
\begin{proof}
By the condition of the lemma, for any $z_0\in E$ and for any
neighborhood $U$ of the point $z_0$ there is a neighborhood
$V\subset U$ such that $(V\setminus E)\cap D^{\,\prime}$ consists of
a finite number of components. Since $x_k\rightarrow z_1$ as
$k\rightarrow\infty,$ there are infinitely many elements of the
sequence $x_k,$ $k=1,2,\ldots ,$ that belong only to one component
of $(V\setminus E)\cap D^{\,\prime}.$ Moreover, all these elements
belong to only one component $K_1$ of the set $D^{\,\prime}\setminus
E,$ because $(V\setminus E)\cap D^{\,\prime}\subset
D^{\,\prime}\setminus E.$ This implies the possibility of choosing
the sequence $x_{k_l}\in K_1,$ $l=1,2,\ldots .$ In order to simplify
notation, we will further assume that the sequence $x_k$ itself
belongs to $K_1.$

\medskip
Let $a_1\in K_1.$ Let us to show the existence of the subsequence
$x_{k_{l_m}},$ $m=1,2,\ldots ,$ and a path $\alpha:[0, 1]\rightarrow
\overline{D^{\,\prime}},$ such that $\alpha(t)\in
D^{\,\prime}\setminus E,$ $x_{k_{l_m}}\in |\alpha|,$ $\alpha(0)=a_1$
and $\alpha(1)=z_1.$ We will use the methodology proposed by
N\"{a}kki and V\"{a}is\"{a}l\"{a} (see \cite[Theorem~17.7
(3)]{Va$_2$} and \cite[Theorem~1.11(3)]{Na$_1$}).

\medskip
Due to the conditions of the lemma, for $U_1=B(z_1, 2^{-1})$ there
is a neighborhood $V_1\subset U_1$ of the point $z_1$ such that
$V_1\cap (D^{\,\prime}\setminus E)$ consists of a finite number of
components. Let us choose among them the one that contains
infinitely many elements of the sequence $x_k,$ $k=1,2,\ldots .$ Let
us denote this component by $W_1$ and fix the element $x_{k_1}$ in
it. Note that, $x_{k_1}\in W_1\cap K_1$ and $W_1$ is connected, so
by definition of the component $K_1$ of the set
$D^{\,\prime}\setminus E$ we have that $W_1\subset K_1.$ Note that,
for $U_2=B(z_1, 2^{-2})$ there is a neighborhood $V_2\subset U_2$ of
$z_1$ such that $V_2\cap W_1$ consists of a finite number of
components. (In the contrary case, for any such $V_2$ the set
$V_2\cap (D^{\,\prime}\setminus E)$ consists of an infinite number
of components, but this contradicts the condition of the lemma). Put
some a component of $V_2\cap W_1,$ which contains an infinite number
of points of the sequence $x_k.$ Let us denote this component by
$W_2,$ a the corresponding point $x_{k_2}\in W_2.$ And so on.
Continuing this process, we will obtain the decreasing sequence of
components
$$K_1\supset
W_1\supset W_2\supset\ldots\,, \qquad W_l\subset B(z_1,
2^{-l})\,,\qquad x_{k_l}\in W_l\,.$$
Each subsequent component is a part of the previous one, so it is
possible in pairs join the points $a_1,$ $x_{k_1},$
$x_{k_2},\ldots,$ $x_{k_l}\,\ldots ,$ by paths $\gamma_l,$
$l=1,2,\ldots ,$ where
\begin{eqnarray*}\gamma_1:[0, 1/2]\rightarrow K_1,\quad \gamma_1(0)=a_1,\quad
\gamma_1(1/2)=x_{k_1} \\
\gamma_2:[1/2, 2/3]\rightarrow W_1,\quad\gamma_2(1/2)=x_{k_1},\quad
\gamma_2(2/3)=x_{k_2}
\end{eqnarray*}
and for $l>2$
$$
\gamma_l:[(l-1)/l, l/(l+1)]\rightarrow W_{l-1},\quad
\gamma_l(l-1/l)=x_{k_{l-1}},\quad\gamma_l(l/(l+1))=x_{k_l}\,.$$
Define the path $\alpha:[0, 1)\rightarrow D^{\,\prime}\setminus E$
as follows:
$$\alpha(t)=\gamma_{l}(t)\,,\quad t\in [(l-1)/l, l/(l+1)]\,.$$
Note that the path $\alpha(t)$ has a continuous extension to a point
$t=1.$ Indeed, let $t_m\rightarrow 1-0.$ Then for any number $l\in
{\Bbb N}$ there is a number $m_l\in {\Bbb N}$ such that
$\frac{l-1}{l}<t_m,$ $m>m_l.$ But then by the construction
$|\gamma_l(t_m)-z_1|<2^{-l+1}.$ We take any $\varepsilon>0.$ Let us
find $l_0=l_0(\varepsilon)\in {\Bbb N}$ such that
$2^{-l_0+1}<\varepsilon.$ Then, due to the mentioned above,
$|\alpha(t_m)-z_1|<2^{-l_0+1}<\varepsilon$ as
$m>m_{l_0}=m_{l_0(\varepsilon)}:=M_0(\varepsilon).$ This proves that
$\alpha(t)\rightarrow z_1$ as $t\rightarrow 1-0.$ Now, $\alpha:[0,
1]\rightarrow D^{\,\prime}\setminus E$ is a desired path.~$\Box$
\end{proof}

\medskip
The key is the following lemma.

\medskip
\begin{lemma}\label{lem2}
{\it\, Let $D^{\,\prime}$ be a domain in ${\Bbb R}^n,$ $n\geqslant
2,$ and let $E\subset \overline{D^{\,\prime}},$ $\partial
D^{\,\prime}\subset E,$ and let $E$ be closed in $\overline{\Bbb
R^{n}}.$ Suppose also that, for of any $z_0\in E$ and for any
neighborhood $U$ of the point $z_0$ there is a neighborhood
$V\subset U$ of this point such that $(V\setminus E)\cap
D^{\,\prime}$ consists of a finite number of components. Let $z_1,
z_2\in E,$ $z_1\ne z_2,$ and let $x_k, y_k\in D^{\,\prime}\setminus
E,$ $k=1,2,\ldots ,$ be sequences, which converge to $z_1$ and
$z_2,$ respectively. Then there are subsequences $x_{k_l}$ and
$y_{k_l},$ $l=1,2,\ldots ,$ each of which belong to exactly one
(corresponding) component $K_1$ and $K_2$ of $D^{\,\prime}\setminus
E.$

Moreover, if $a_1\in K_1$ and $a_2\in K_2,$ $a_1\ne a_2,$ then there
are (possibly) new subsequences $x_{k_{l_m}}$ and $y_{k_{l_m}},$
$m=1,2,\ldots ,$ and paths $\alpha:[0, 1]\rightarrow
\overline{D^{\,\prime}},$ $\beta:[0, 1]\rightarrow
\overline{D^{\,\prime}}$ such that:

\medskip
a) $\alpha(t)\in D^{\,\prime}\setminus E$ and $\beta(t)\in
D^{\,\prime}\setminus E$ for $0\leqslant t<1,$

\medskip
b) $x_{k_{l_m}}\in |\alpha|,$ $y_{k_{l_m}}\in |\beta|,$
$m=1,2,\ldots ,$

\medskip
c) $\alpha(0)=a_1,$ $\beta(0)=a_2,$ $\alpha(1)=z_1,$ $\beta(1)=z_2,$

\medskip
d) $|\alpha|\cap |\beta|=\varnothing,$ see Figure~\ref{fig3}.
\begin{figure}[h]
\centerline{\includegraphics[scale=0.4]{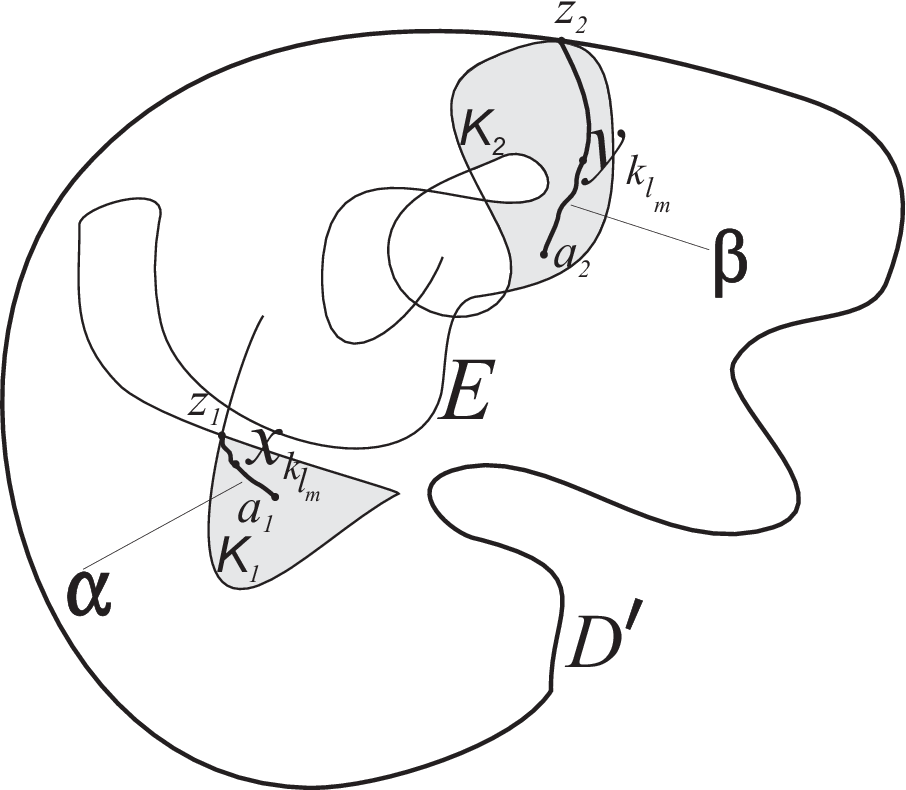}} \caption{The
formulation of Lemma~\ref{lem2}}\label{fig3}
\end{figure}}
\end{lemma}

\medskip
\begin{proof}
The existence of paths $\alpha$ and $\beta,$ satisfying conditions
a)--c), follows by Lemma~\ref{lem4}. It remains only to check the
condition $d),$ namely, to explain why the paths $\alpha$ and
$\beta$ may be chosen disjoint. Two options are possible:

\medskip
1) $K_1\ne K_2,$ i.e., paths $\alpha$ and $\beta$ (except of their
right endpoints) belong to (disjoint) components of
$D^{\,\prime}\setminus E.$ Then $\alpha$ and $\beta$ are disjoin by
the definition.

\medskip
2) $K_1=K_2.$ Then the paths $\alpha$ and $\beta$ excluding their
right endpoints belong to the domain $D:=K_1=K_2,$ for which points
$z_1$ and $z_2$ are boundary points. Let us discard no more than a
countable number of loops in the path $\alpha$ and denote the
corresponding new path by $\alpha_0.$ In this case, the possibility
of choosing disjoint paths $\alpha^{\,*}$ and $\beta^{\,*},$ which
satisfy conditions a), c) and d) is a statement of Lemma~\ref{lem3},
and due to this lemma it is possible set $\alpha^{\,*}=\alpha_0,$
and $\beta^{\,*}(t)$ coincides with $\beta(t)$ for $t>T_0$ and some
$0<T_0<1.$ It remains to explain, why it is possible to choose paths
such that condition b) is also fulfilled, i.e., $x_{k_{l_m}}\in
|\alpha^{\,*}|,$ $y_{k_{l_m}}\in |\beta^{\,*}|,$ $m=1,2,\ldots .$
Regarding $\beta^{\,*},$ the latter is obvious, because
$\beta^{\,*}(t)$ coincides with $\beta(t)$ for $t>T_0$ and some
$0<T_0<1.$ In particular, this means that all elements $y_{k_{l_m}}$
belong to $\beta,$ starting from some number, $m\geqslant M_0.$

In addition, since $\beta^{\,*}$ bypasses the point $z_1,$ there
exists a neighborhood $U$ of $z_1$ such that $U\cap
|\beta^{\,*}|=\varnothing.$ On the other hand, since
$\alpha\rightarrow z_1$ as $t\rightarrow 1,$ there exists $0<t_1<1$
such that $\alpha(t)\in U$ at $t>t_1.$ Let
$M_1:=\min\limits_{m\in {\Bbb N}: x_{k_{l_m}}\in
|\alpha(t)|_{[t>t_1]|}}m.$
Then we add to the path $\alpha_0$ all discarded loops of the path
$\alpha,$ which correspond to the elements $x_{k_{l_m}}$ for
$m\geqslant M_1.$ We denote the last path by $\alpha^{\,**}.$

\medskip
The paths $\alpha^{\,**}$ and $\beta^{\,*}$ are required: they
satisfy the conditions a), c) and d), in addition, each of these
paths contains itself, respectively, all elements $x_{k_{l_m}}$ and
$y_{k_{l_m}}$ for $m\geqslant \max\{M_0, M_1\}.$ (We may go to
subsequences so that $m$ changes starting from~1). Lemma is
proved.~$\Box$
\end{proof}

\medskip
{\it Proof of Theorem~\ref{th2}.} {\bf I.} Let $f\in {\frak S}^P_{E,
\delta, Q}(D, D^{\,\prime}).$ By Theorem~\ref{th3} the mapping $f$
have a continuous extension $\overline{f}:\overline{D}\rightarrow
\overline{D^{\,\prime}}$ such that
$\overline{f}(\overline{D})=\overline{D^{\,\prime}}.$ The
equicontinuity of the family of mappings ${\frak S}^P_{\delta, A, Q
}(\overline{D}, \overline{D^{\,\prime}})$ in $D$ is a statement of
Theorem~1.1 in \cite{SevSkv$_4$}. It remains to establish the
equicontinuity of this family in $\partial D.$

We will carry out the proof by the opposite. Assume that there is
$x_0\in \partial D,$ a number $\varepsilon_0>0,$ a sequence $x_m\in
\overline{D},$ which converges to the point $x_0$ as
$m\rightarrow\infty,$ and the corresponding maps $\overline{f}_m\in
{\frak S}^ P_{E, \delta, Q }(\overline{D}, \overline{D})$ such that
\begin{equation}\label{eq12A}
h(\overline{f}_m(x_m),\overline{f}_m(x_0))\geqslant\varepsilon_0,\quad
m=1,2,\ldots .
\end{equation}
Let us put $f_m:=\overline{f}_m|_{D}.$ Since $f_m$ has a continuous
extension to $\partial D,$ we may assume that $x_m\in D.$ Therefore,
$\overline{f }_m(x_m)=f_m(x_m).$ In addition, there is a sequence
$y_m\in D$ such that $y_m\rightarrow x_0$ as $m\rightarrow\infty$
and $h(f_m(y_m),\overline{f}_m(x_0))\rightarrow 0$ as
$m\rightarrow\infty.$ Since the space $\overline{{\Bbb R}^n}$ is
compact, we may assume that $f_m(x_m)$ and $\overline{f}_m(x_0)$
converge as $m\rightarrow\infty$ to some points $z_1$ and $z_2.$ By
the continuity of the chordal metric $h,$ (\ref{eq12A}) implies that
$z_1\ne z_2.$ Note that, $z_2\in E,$ since $\overline{f}_m(x_0)\in
C(f_m, x_0)\subset E$ and $E$ is closed by the definition of the
class ${\frak S}^P_{E, \delta, Q}(D, D^{\,\prime}).$ However, we do
not know if $z_1$ is in $E.$

\medskip
{\bf II.} Let us show that we may also choose $x_m, y_m\in D$ so
that $f_m(x_m), f_m(y_m)\not\in E$ for all $m\in {\Bbb N}.$ Indeed,
since under condition~2) the set $f_m^{\,-1}(E\cap D^{\,\prime})$ is
nowhere dense in $D,$ there is a sequence $x_{km}\in
D\setminus(f^{\,-1}(E)\cap D^{\,\prime}),$ $k=1,2,\ldots ,$ such
that $x_{km}\rightarrow x_m$ as $k\rightarrow\infty.$ Fix
$\varepsilon>0.$ Due to the continuity of the mapping $f_m$ at the
point $x_m,$ given a number $m\in {\Bbb N}$ there exists a number
$k_m\in {\Bbb N}$ such that $h(f_m(x_{k_mm}),
f_m(x_m))<\frac{1}{2^m}.$ Therefore, by the triangle inequality
$$h(f_m(x_{k_mm}), z_1)\leqslant h(f_m(x_{k_mm}), f_m(x_m))+
h(f_m(x_m), z_1)\leqslant \frac{1}{2^m}+\varepsilon\,,$$
$m\geqslant m_0=m_0(\varepsilon),$ because $f_m(x_m)\rightarrow z_1$
as $m\rightarrow\infty$ by the choice of $x_m$ and $z_1.$ Therefore,
$x_m\in D$ may be replaced by $x_{k_mm}\in
D\setminus(f_m^{\,-1}(E\cap D^{\,\prime}))$ as needed. We may reason
similarly for~$y_m.$

\medskip
{\bf III.} Now, we construct some paths $\alpha$ and $\beta$ that
contain sequences $f_m(x_m)$ and $f_m(y_m),$ some points $a_1$ and
$a_2\in P,$ and have positive distance from each other. Let us
consider two cases:

\medskip
{\bf III. a)} If $z_1\in E.$ Due to Lemma~\ref{lem2} it is possible
to chose subsequences $f_{m_l}(x_{m_l})$ and $f_{m_l}(y_{m_l}),$
$l=1,2,\ldots ,$ each of which belongs to exactly one
(corresponding) component $K_1$ and $K_2$ of the set
$D^{\,\prime}\setminus E.$ Let $a_{i_1}:=x_{K_1}\in K_1\cap P$ and
$a_{i_2}:=y_{K_2}\in K_2\cap P$ be two different points of the set
$P$ (where the situation $K_1=K_2$ is not excluded) such that
$$h(f_{m_l}^{-1}(x_{K_1}),
\partial D)=h(f_{m_l}^{-1}(a_{i_1}),
\partial D)\geqslant \delta$$
and
$$h(f_{m_l}^{-1}(y_{K_2}),
\partial D)=h(f_{m_l}^{-1}(a_{i_2}),
\partial D)\geqslant \delta\,.$$
The existence of such points follows by the definition of the class
${\frak S}^P_{E, \delta, Q}(D, D^{\,\prime}).$ By Lemma~\ref{lem2}
there are new subsequences $f_{m_{l_k}}(x_{m_{l_k}})$ and
$f_{m_{l_k}}(y_{m_{l_k}}),$ $k=1,2,\ldots ,$ and paths $\alpha:[0,
1]\rightarrow \overline{D^{\,\prime}},$ $\beta:[0, 1]\rightarrow
\overline{D^{\,\prime}}$ such that:

\medskip
a) $\alpha(t)\in D^{\,\prime}\setminus E$ and $\beta(t)\in
D^{\,\prime}\setminus E$ for $0\leqslant t<1,$

\medskip
b) $f_{m_{l_k}}(x_{m_{l_k}})\in |\alpha|,$
$f_{m_{l_k}}(y_{m_{l_k}})\in |\beta|,$ $k=1,2,\ldots ,$

\medskip
c) $\alpha(0)=a_{i_1},$ $\beta(0)=a_{i_2},$ $\alpha(1)=z_1,$
$\beta(1)=z_2,$

\medskip
d) $|\alpha|\cap |\beta|=\varnothing.$

\begin{figure}[h]
\centerline{\includegraphics[scale=0.5]{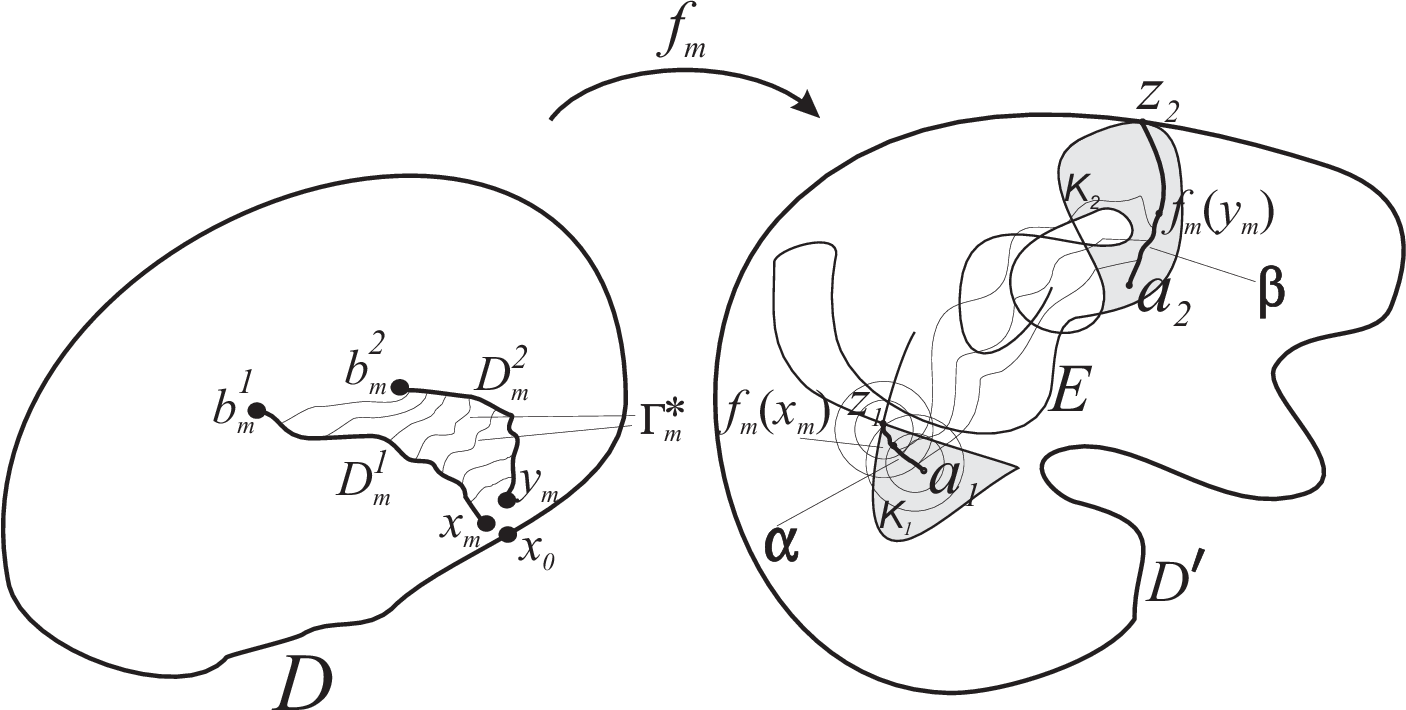}} \caption{To
the proof of Theorem~\ref{th2}, case {\bf III. a)}}\label{fig5}
\end{figure}

Without loss of generality, we may assume that the sequences
$f_m(x_m)$ and $f_m(y_m)$ satisfy the above conditions a)--d), in
addition, we may assume that $a_{i_1}=a_1$ and $a_{i_2}=a_2,$ see
Figure~\ref{fig5}.
Let $$f_m(x_m)=\alpha(t_m)\,, \qquad f_m(y_m)=\beta(p_m)\,, \quad
0<t_m, p_m<1\,.$$
{\bf III. b)} If $z_1\not\in E.$ Since by the assumption $\partial
D^{\,\prime}\subset E,$ the point $z_1$ belongs to some component
$K_1$ of the set $D^{\,\prime}\setminus E.$ Since $E$ is closed in
$\overline{{\Bbb R}^n},$ it is closed with respect to
$D^{\,\prime}.$ Therefore, $K_1$ is an open set. Then the sequence
$f_m(x_m),$ which converges to $z_1$ as $m\rightarrow\infty,$ itself
belongs to $K_1$ starting from some number, $m>M_1.$ Let $a_{i_1}\in
K_1\cap P$ be a point from the condition of the theorem. Without
loss of generality, we may assume that $a_{i_1}=a_1.$ Consequently
joining the points $f_{M_1}(x_{M_1}),$ $f_{M_1+1}(x_{M_1+1}),$
$\ldots,$ $f_{m}(x_{m}),$ $\ldots,$ going to the subsequence if it
needs and reasoning as under the proof of Lemma~\ref{lem4}, we
obtain a path $\widetilde{\alpha}:[1/2, 1]\rightarrow K_1$ starting
at the point $f_{M_1}(x_{M_1})$ and ending at the point $z_1,$
$\widetilde{\alpha}(t_m)=f_m(x_m),$ $t_m\in [1/2, 1],$
$m=1,2,\ldots, $ $t_m\rightarrow 1$ as $m\rightarrow \infty,$ which
lies completely in some ball $B(z_1, \varepsilon_1)\subset
D^{\,\prime}$ and converges to $z_1$ as $t\rightarrow 1-0.$

\medskip
By Lemma~\ref{lem4} there exists a subsequence $f_{m_l}(y_{m_l}),$
$l=1,2,\ldots ,$ which belongs to exactly one component $K_2$ of the
set $D^{\,\prime}\setminus E.$ Let $a_{i_2}\in K_2\cap P$ be a point
from the condition of the class ${\frak S}^P_{E, \delta, Q}(D,
D^{\,\prime}).$ Then there is a new subsequence
$f_{m_{l_k}}(y_{m_{l_k}}),$ $k=1,2,\ldots ,$ and a path $\beta:[0,
1]\rightarrow \overline{D^{\,\prime}}$ such that:

\medskip
a) $\beta(t)\in D^{\,\prime}\setminus E$ for $0\leqslant t<1,$

\medskip
b) $f_{m_{l_k}}(y_{m_{l_k}})\in |\beta|,$ $m=1,2,\ldots ,$

\medskip
c) $\beta(0)=a_{i_2},$ $\beta(1)=z_2.$

\medskip
To simplify notation, we believe that these conditions are satisfied
by itself sequence $f_m(y_m)$ instead of $f_{m_{l_k}}(y_{m_{l_k}}),$
moreover, $a_{i_2}=a_2.$

\begin{figure}[h]
\centerline{\includegraphics[scale=0.5]{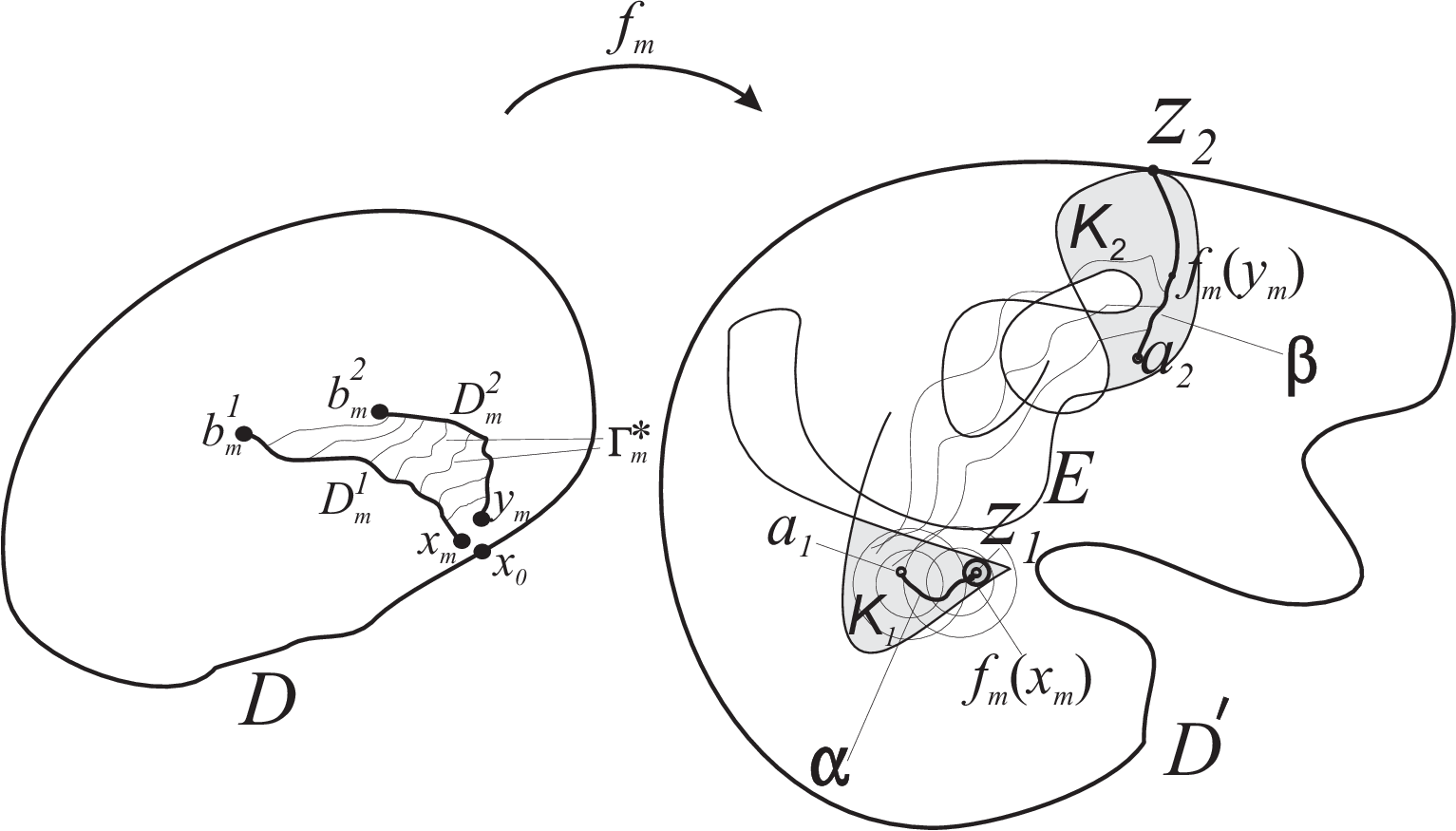}} \caption{To
the proof of Theorem~\ref{th2}, case {\bf III. b)}}\label{fig6}
\end{figure}

If $z_1\ne a_1,$ join the point $a_1$ with the point
$f_{M_1}(x_{M_1})$ by a path $\Delta:[0, 1/2]\rightarrow K_1$ and
put
$$\alpha(t):=\begin{cases}\Delta(t)\,,& t\in [0, 1/2]\,,\\
\widetilde{\alpha}(t)\,,& t\in [1/2, 1]\,.\end{cases}$$
If $z_1=a_1,$ the procedure is similar, but we choose $\Delta$ so
that it belongs entirely to the ball $B(z_1, \varepsilon_1).$ In
addition, we may choose the radius $\varepsilon_1$ of the ball
$B(z_1, \varepsilon_1)$ so small that $|\beta|\cap B(z_1,
\varepsilon_1)=\varnothing,$
see Figure~\ref{fig6}.
Therefore, reasoning as in the proof of Lemma~\ref{lem2}, it may be
shown that the paths $\alpha$ and $\beta$ may be chosen disjoint (if
$z_1=a_1,$ this follows from the construction of these paths. If
$z_1\ne a_1,$ we first move from $\alpha$ its loops and apply
Lemma~\ref{lem3}. As a result we get a pair of disjoint paths, one
of them is Jordan curve, and another contains <<start>> and <<end>>
of $\beta,$ whenever $\beta$ already contains the corresponding
points $f_m(y_m).$ Finally, we <<return>> to $\alpha$ some part of
its loops. The resulting path also contains $f_m(x_m)$ for
sufficiently large $m\in {\Bbb N}$ and does not intersect with
$\beta,$ as well).

\medskip
{\bf IV.} Put
$$C^1_m(t)=\alpha(1-t)|_{[1-t_m, 1]}\,,\qquad C^2_m(t)=\beta(1-t)|_{[ 1-p_m, 1]}\,.$$
Without loss of the generality, we may assume that $z_1\ne\infty.$
Otherwise, we relabel $z_1$ and $z_2.$ Let $D^1_m: [1-t_m,
c^1_m)\rightarrow D$ and $D^2_m: [1-p_m, c^2_m)\rightarrow D$ be
maximal $f_m$-liftings of paths $C^1_m$ and $C^2_m$ starting at
points $x_m$ and $y_m$ (they exist by Proposition~\ref{pr3}). By the
same proposition, one of two situations is possible:

\medskip
1) $D^1_m(t)\rightarrow x_1\in D$ as $t\rightarrow c^1_m-0$ and
$c^1_m=1$ and $f_m(D^1_m(1))=f_m(x_1)=a_1,$

\medskip
2) or $D^1_m(t)\rightarrow \partial D$ as $t\rightarrow c^1_m.$

\medskip
Let us show that situation~2) does not hold. Assume the opposite:
let $D^1_m(t)\rightarrow \partial D$ as $t\rightarrow c^1_m.$ We
choose an arbitrary sequence $\theta_k\in [1-t_m, c^1_m)$ such that
$\theta_k\rightarrow c^1_m-0$ as $k\rightarrow\infty.$ Since the
space $\overline{{\Bbb R}^n}$ is compact, $\partial D$ is also
compact as a closed subset of a compact space. Then there exists
$w_k\in
\partial D$ is such that
\begin{equation}\label{eq7C}
h(D^1_m(\theta_k), \partial D)=h(D^1_m(\theta_k), w_k) \rightarrow
0\,,\qquad k\rightarrow \infty\,.
\end{equation}
Due to the compactness of $\partial D$, we may assume that
$w_k\rightarrow w_0\in \partial D$ as $k\rightarrow\infty.$
Therefore, from the relation~(\ref{eq7C}) and by the triangle
inequality
\begin{equation}\label{eq8C}
h(D^1_m(\theta_k), w_0)\leqslant h(D^1_m(\theta_k), w_k)+h(w_k,
w_0)\rightarrow 0\,,\qquad k\rightarrow \infty\,.
\end{equation}
On the other hand,
\begin{equation}\label{eq9C}
f_m(D^1_m(\theta_k))=C^1_m(\theta_k)=\alpha(1-\theta_k)\rightarrow
\alpha(1-c^1_m) \,,\quad k\rightarrow\infty\,,
\end{equation}
because by the construction the whole path
$C^1_m=\alpha(1-t)|_{[1-t_m, 1]}$ lies in $D^{\,\prime}\setminus
E\subset D^{\,\prime}\setminus C(f_m,\partial D)$ together with its
endpoints. At the same time, by~(\ref{eq8C}) and~(\ref{eq9C}) we
have that $\alpha(1-c^1_m)\in C(f_m, \partial D).$ The inclusion
$|\alpha(1-t)|_{[1-t_m, c^1_m]}|\subset D^{\,\prime}\setminus C(f_m,
\partial D)$ contradicts the relation $\alpha(1-c^1_m)\in C(f_m,
\partial D).$ The resulting contradiction indicates the
impossibility of the second situation above.

\medskip
Therefore, the first situation is fulfilled: $D^1_m(t)\rightarrow
x_1\in D$ as $t\rightarrow c^1_m-0$ and $c^1_m=1$ and $f_m(D^1_m(1))
=f_m(x_1)=a_1.$ In other words, the lifting $D^1_m$ is complete,
i.e., $D^1_m:[1-t_m, 1]\rightarrow D.$ Similarly, the path $C^2_m$
has a complete $f_m$-lifting $D^2_m: [1-p_m, 1]\rightarrow D.$ By
the definition of ${\frak S}^P_{E, \delta, Q}(D, D^{\,\prime}),$ we
have that $h(f_m^{\,-1}(a_i),
\partial D)\geqslant \delta,$ $i=1,2 .$ Thus, the end points of
$D^1_m$ and $D^2_m,$ which we will denote by $b_m^1$ and $b_m^2,$
are distant from the boundary of $D$ not less than on $\delta.$

As always, we denote by $|C^1_m|$ and $|C^2_m|$ the loci of $C^1_m$
and $C^2_m,$ respectively. Put
$l_0={\rm dist}\,(|\alpha|, |\beta|)$
and consider the coverage $A_0:=\bigcup\limits_{x\in |\alpha|}B(x,
l_0/4)$ of $|\alpha|$ using balls. Since $|\alpha|$ is a compact
set, we may choose a finite number of indices $1\leqslant
N_0<\infty$ and corresponding points $\omega^{\,*}_1,\ldots,
\omega^{\,*}_{N_0}\in |\alpha|$ such that $|\alpha|\subset
B_0:=\bigcup\limits_{i=1}^{N_0}B(\omega^{\,*}_i, l_0/4).$
In this case,
$|C^1_m|\subset \bigcup\limits_{i=1}^{N_0}B(\omega^{\,*}_i, l_0/4).$
Let $\Gamma_m$ be a family of all paths joining $|C^1_m|$ and
$|C^2_m|$ in $D^{\,\prime}.$ Now, we have
\begin{equation}\label{eq10C}
\Gamma_m=\bigcup\limits_{i=0}^{N_0}\Gamma_{mi}\,,
\end{equation}
where $\Gamma_{mi}$ is a family of all paths $\gamma:[0,
1]\rightarrow D^{\,\prime}$ such that $\gamma(0)\in
B(\omega^{\,*}_i, l_0/4)\cap |C^1_m|$ and $\gamma(1)\in |C_2^m|$ for
$1\leqslant i\leqslant N_0.$ Due to~\cite[Theorem~1.I.5.46]{Ku}, we
may show that
\begin{equation}\label{eq11C}
\Gamma_{mi}>\Gamma(S(\omega^{\,*}_i, l_0/4), S(\omega^{\,*}_i,
l_0/2), A(\omega^{\,*}_i, l_0/4, l_0/2))\,.
\end{equation}
Let $$\Gamma^{\,*}_m:=\Gamma(|D_m^1|, |D_m^2|, D)\,.$$
Note that $f_m(\Gamma^{\,*}_m)\subset\Gamma_m.$ Then, by
(\ref{eq10C}) and (\ref{eq11C})
\begin{equation}\label{eq6A}
\Gamma^{\,*}_m>\bigcup\limits_{i=1}^{N_0}\Gamma_{f_m}(\omega^{\,*}_i,
l_0/4, l_0/2)\,.
\end{equation}%
Set $\widetilde{Q}(y)=\max\{Q(y), 1\},$ and let
$$\widetilde{q}_{\omega^{\,*}_i}(r)=\frac{1}{\omega_{n-1}r^{n-1}}
\int\limits_{S(\omega^{\,*}_i,
r)}\widetilde{Q}(y)\,d\mathcal{H}^{n-1}\,.$$
By the assumption, $q_{\omega^{\,*}_i}(r)\ne \infty$ for $r\in E_1,$
where $E_1$ is some set of positive measure in $[l_0/4, l_0/2].$
Thus, $\widetilde{q}_{\omega^{\,*}_i}(r)\ne \infty$ for $r\in E_1,$
as well.
Let
\begin{equation}\label{eq2}
I_i=I_i(\omega^{\,*}_i, l_0/4, l_0/2)=\int\limits_{l_0/4}^{l_0/2}\
\frac{dr}{r\widetilde{q}_{\omega^{\,*}_i}^{\frac{1}{n-1}}(r)}\,.
\end{equation}
Note that, $I_i\ne 0,$ because there is some strictly positive
function (on some set of positive measure) in the integral
in~(\ref{eq2}). In addition, $I_i\ne\infty,$ because
$I_i\leqslant\log\frac{r_2}{r_1}<\infty,$ $i=1,2, \ldots, {N_0}.$
Let us put now
$$\eta_i(r)=\begin{cases}
\frac{1}{I_ir\widetilde{q}_{\omega^{\,*}_i}^{\frac{1}{n-1}}(r)}\,,&
r\in [l_0/4, l_0/2]\,,\\
0,& r\not\in [l_0/4, l_0/2]\,.
\end{cases}$$
Note that the function ~$\eta_i$ satisfies the condition
~$\int\limits_{l_0/4}^{l_0/2}\eta_i(r)\,dr=1,$ therefore, it may be
substituted in the right-hand side of~(\ref{eq2*A}). Now, we obtain
that
\begin{equation}\label{eq7BA}
\Gamma_{f_m}(\omega^{\,*}_i, l_0/4, l_0/2)\leqslant
\int\limits_{A(\omega^{\,*}_i, l_0/4, l_0/2)}
\widetilde{Q}(y)\,\eta^n_i(|y-\omega^{\,*}_i|)\,dm(y)\,.\end{equation}
However, by the Fubini theorem (see, e.g.,
\cite[theorem~8.1.III]{Sa})
\begin{multline}\label{eq7CB}
\int\limits_{A(\omega^{\,*}_i, l_0/4, l_0/2)}
\widetilde{Q}(y)\,\eta^n_i(|y-\omega^{\,*}_i|)\,dm(y)
=\int\limits_{l_0/4}^{l_0/2}\int\limits_{S(\omega^{\,*}_i,
r)}\widetilde{Q}(y)\eta^n_i(|y-\omega^{\,*}_i|)\,d\mathcal{H}^{n-1}\,dr\,
\\=\frac{\omega_{n-1}}{I_i^n}\int\limits_{l_0/4}^{l_0/2}r^{n-1}
\widetilde{q}_{\omega^{\,*}_i}(r)\cdot
\frac{dr}{r^n\widetilde{q}^{\frac{n}{n-1}}_{\omega^{\,*}_i}(r)}=
\frac{\omega_{n-1}}{I_i^{n-1}}\,.\end{multline}
Then from~(\ref{eq7BA}) and~(\ref{eq7CB}) it follows that
$$M(\Gamma_{f_m}(\omega^{\,*}_i, l_0/4, l_0/2))\leqslant \frac{\omega_{n-1}}{I_i^{n-1}}\,,$$
whence from~(\ref{eq6A})
\begin{equation}\label{eq14A}
M(\Gamma^{\,*}_m)\leqslant
\sum\limits_{i=1}^{N_0}M(\Gamma_{f_m}(\omega^{\,*}_i, l_0/4,
l_0/2))\leqslant
\sum\limits_{i=1}^{N_0}\frac{\omega_{n-1}}{I_i^{n-1}}:=C_0\,, \quad
m=1,2,\ldots\,.
\end{equation}
We will show that the relation~(\ref{eq14A}) contradicts the
condition of weak flatness of $\partial D.$ Indeed, by the
construction
$$h(|D^1_m|)\geqslant h(x_m, b_m^1) \geqslant
(1/2)\cdot h(b_m^1, \partial D)>\delta/2\,,$$
\begin{equation}\label{eq15}
h(|D^2_m|)\geqslant h(y_m, b_m^2) \geqslant (1/2)\cdot h(b_m^2,
\partial D)>\delta/2
\end{equation}
for all $m\geqslant M_0$ and for some $M_0\in {\Bbb N}.$
Let us put $U:=B_h(x_0, r_0)=\{y\in \overline{{\Bbb R}^n}: h(y,
x_0)<r_0\},$ where $0<r_0<\delta/4$ and a number $\delta$ refers to
relation~(\ref{eq15}). Note that, $|D^1_m|\cap U\ne\varnothing\ne
|D^1_m|\cap (D\setminus U)$ for each $m\in{\Bbb N},$ because
$h(|D^1_m|)\geqslant \delta/2$ and $x_m\in |D^1_m|,$ $x_m\rightarrow
x_0$ as $m\rightarrow\infty.$ Similarly, $|D^2_m|\cap
U\ne\varnothing\ne |D^2_m|\cap (D\setminus U).$ Since $|D^1_m|$ and
$|D^2_m|$ are continua,
\begin{equation}\label{eq8A}
|D^1_m|\cap \partial U\ne\varnothing, \quad |D^2_m|\cap
\partial U\ne\varnothing\,,
\end{equation}
see, e.g.,~\cite[Theorem~1.I.5.46]{Ku}. Since $\partial D$ is weakly
flat, for $P:=C_0>0$ (where $C_0$ is a number in~(\ref{eq14A}))
there is a neighborhood $V\subset U$ of $x_0$ such that
\begin{equation}\label{eq9A}
M(\Gamma(E, F, D))>C_0
\end{equation}
for any continua $E, F\subset D$ such that $E\cap
\partial U\ne\varnothing\ne E\cap \partial V$ and $F\cap \partial
U\ne\varnothing\ne F\cap \partial V.$ Let us show that,
\begin{equation}\label{eq10A}
|D^1_m|\cap \partial V\ne\varnothing, \quad |D^2_m|\cap
\partial V\ne\varnothing\end{equation}
for sufficiently large $m\in {\Bbb N}.$ Indeed, $x_m\in |D^1_m|$ and
$y_m\in |D^2_m|,$ where $x_m, y_m\rightarrow x_0\in V$ as
$m\rightarrow\infty.$ In this case, $|D^1_m|\cap V\ne\varnothing\ne
|D^2_m|\cap V$ for sufficiently large $m\in {\Bbb N}.$ Note that,
$h(V)\leqslant h(U)\leqslant 2r_0<\delta/2.$ By~(\ref{eq15}),
$h(|D^1_m|)>\delta/2.$ Therefore, $|D^1_m|\cap (D\setminus
V)\ne\varnothing$ and, therefore, $|D^1_m|\cap\partial
V\ne\varnothing$ (see, e.g.,~\cite[Theorem~1.I.5.46]{Ku}).
Similarly, $h(V)\leqslant h(U)\leqslant 2r_0<\delta/2.$
From~(\ref{eq15}) it follows that $h(|D^2_m|)>\delta/2,$ therefore,
$|D^2_m|\cap (D\setminus V)\ne\varnothing.$
By~\cite[Theorem~1.I.5.46]{Ku} we obtain that $|D^2_m|\cap\partial
V\ne\varnothing.$ Thus, the relation~(\ref{eq10A}) is established.
Combining relations~(\ref{eq8A}), (\ref{eq9A}) and (\ref{eq10A}), we
obtain that $M(\Gamma^{\,*}_m)=M(\Gamma(|D^1_m|, |D^2_m|, D))>C_0.$
The last condition contradicts~(\ref{eq14A}), which proves the
theorem.~$\Box$

\medskip
{\bf Funding.} The work was supported by the National Research
Foundation of Ukraine (Project ``Analogues of Carath\'{e}odory and
Koebe-Bloch theorems for Orlycz-Sobolev classes'', Project number
2025.02/0010).

\medskip
{\bf \noindent Victoria Desyatka} \\
Zhytomyr Ivan Franko State University,  \\
40 Velyka Berdychivs'ka Str., 10 008  Zhytomyr, UKRAINE \\
victoriazehrer@gmail.com

\medskip
\medskip
{\bf \noindent Evgeny Sevost'yanov} \\
{\bf 1.} Zhytomyr Ivan Franko State University,  \\
40 Velyka Berdychivs'ka Str., 10 008  Zhytomyr, UKRAINE \\
{\bf 2.} Institute of Applied Mathematics and Mechanics\\
of NAS of Ukraine, \\
19 Henerala Batyuka Str., 84 116 Slov'yansk,  UKRAINE\\
esevostyanov2009@gmail.com

\end{document}